\documentclass[11pt]{amsart}
\usepackage {amsmath, amssymb, a4wide, epsfig, enumerate, psfrag}
\usepackage[latin1]{inputenc}
\usepackage[english]{babel}

\date{\today}

 \author{Romain Dujardin}

\address{CMLS \\ \'Ecole Polytechnique \\ 91128 Palaiseau\\
         France}
\email{dujardin@math.polytechnique.fr}
\thanks{Research partially supported by ANR project BERKO}

\title[Fatou directions]{Fatou directions along the Julia set for endomorphisms of $\mathbb{CP}^k$}

\newcommand{\cc}{\mathbb{C}}

\newcommand{\re}{\mathbb{R}}

\newcommand{\nn}{\mathbb{N}}
\newcommand{\pp}{\mathbb{P}}
\newcommand{\e}{\varepsilon}
\newcommand{\cv}{\rightarrow}

\newcommand{\fr}{\partial}
\newcommand{\om}{\Omega}
\newcommand{\set}[1]{\left\{#1\right\}}
\newcommand{\norm}[1]{\left\Vert#1\right\Vert}
\newcommand{\abs}[1]{\left\vert#1\right\vert}
\newcommand{\cd}{{\cc^2}}
\newcommand{\pd}{{\mathbb{P}^2}}
\newcommand{\pu}{{\mathbb{P}^1}}
\newcommand{\pk}{{\mathbb{P}^k}}
\newcommand{\rest}[1]{ \arrowvert_{#1}}
\newcommand{\m}{{\bf M}}
\newcommand{\unsur}[1]{\frac{1}{#1}}

\newcommand{\qq}{\mathcal{Q}}

\renewcommand{\o}{\omega}
\newcommand{\lrpar}[1]{\left(#1\right)}
\newcommand{\bra}[1]{\left\langle #1\right\rangle}
\newcommand{\F}{\mathcal{F}}
\newcommand{\cC}{\mathcal{C}}
\newcommand{\la}{\lambda}

\newcommand{\itm}{\item[-]}

\DeclareMathOperator{\supp}{Supp}
\DeclareMathOperator{\vol}{Vol}

\DeclareMathOperator{\dist}{dist}

\DeclareMathOperator{\rank}{rank}
\DeclareMathOperator{\trace}{trace}
\DeclareMathOperator{\leb}{Leb}

\newtheorem{prop} {Proposition} [section]
\newtheorem{thm}[prop] {Theorem} 
\newtheorem{defi}[prop] {Definition}
\newtheorem{lem}[prop] {Lemma}
\newtheorem{cor}[prop]{Corollary}
\newtheorem{conj}[prop]{Conjecture}

\newtheorem{question}[prop]{Question}
\theoremstyle{remark}

\newtheorem{rmk}[prop]{Remark}

\begin{document}

\begin{abstract} We study the dynamics on the Julia set for holomorphic endomorphisms of $\mathbb{CP}^k$. The Julia set is the suppport of the so-called Green current $T$, so it admits a  natural filtration
$J=J_1\supset \cdots\supset J_k$, where for $1\leq q\leq k$ we put $J_q = \mathrm{Supp}(T^q)$.
 We show that for a generic point of $J_q\setminus J_{q+1}$ there are at least $(k-q)$ ``Fatou directions" in the tangent space. We also give estimates for the rate of expansion in directions transverse to the Fatou directions.
\end{abstract}

\maketitle

\section*{Introduction}

\subsection{Background}
In this paper we are concerned with iteration theory of holomorphic endomorphisms
of the complex projective space in several dimensions. To fix notation, let $k>1$ and $f: \pk \cv\pk$ be such an endomorphism, and  $d\geq 2$ be its degree. Recall that $d$ is the degree of the hypersurface $f^{-1}(H)$, where $H$ is a generic hyperplane. The topological degree of $f$ is $d^k$.

A great  achievement in this area of research is
the construction and study of the so-called {\em equilibrium measure} $\mu$, in particular through the work of
Hubbard, Papadopol, Forn\ae ss, Sibony, Briend, Duval and Dinh \cite{hp, fs2,briend-duval1,briend-duval2,ds-pl}. In dimension 1, these results were previously  obtained independently by Lyubich \cite{ly} and Freire-Lopes-Ma\~né \cite{flm}.
The equilibrium measure is defined as the limit of the sequence of measures $\mu_n =\unsur{d^{kn}} \sum_{ f^n(y) = x} \delta_y$, where $x\in \pk$ is a generic point. Among many interesting dynamical properties, let us only mention that $\mu$ is {\em repelling}, in the sense that its Lyapunov exponents     are greater than or equal to $\frac{\log d}{2}$ and that it describes the asymptotic distribution of repelling periodic orbits   \cite{briend-duval1}.

\medskip

On the other hand, the basic understanding of the dynamics outside $\supp(\mu)$ remains problematic.
Let us classically  define the Fatou set $F$ as the (open)
set of points where $(f^n)_{n\geq 0}$  locally defines  a normal family of mappings, and the Julia set
by $J = \pk\setminus F$. As opposed to dimension 1, the Julia set
 is usually larger than  $\supp(\mu)$. Our main purpose in this article is to study the structure of the  dynamics on $J\setminus \supp(\mu)$.

The Julia set is the support of another invariant measurable object: the {\em Green current} $T$
\cite{hp, fs2}. It is a closed positive current of bidegree $(1,1)$
 defined as follows:  if $H$ is a generic hyperplane,
 $T = \lim_{n\cv \infty} d^{-n} [f^{-n}(H)]$. Here, as usual, the notation $[\cdot]$ stands for the integration current.  The fact that $J=\supp(T)$ was proven independently by Forn\ae ss and Sibony \cite{fs2} and Ueda \cite{ueda-fatou}. Furthermore, $T$ has continuous potential, so its exterior powers are  well-defined and, in fact, $\mu=T^k$.

\medskip

For simplicity, let us temporarily work in dimension $k=2$. Put $J_1=J$ and $J_2 = \supp \mu$ so that $J_2\subset J_1$. We say that a holomorphic disk $\Delta\subset\pd$  is a {\em Fatou disk} if $(f^n\rest{\Delta})_{n\geq 0}$ is a normal family.
Many authors have suggested  to understand the difference between the dynamics on $J_2$ and $J_1\setminus J_2$
by the presence of Fatou disks  ``filling''  $J_1\setminus J_2$ (in an appropriate sense). This issue already appears in \cite{fs2}.
Thus the dynamics on $J_1\setminus J_2$ would be in a sense Fatou in the ``tangential'' direction and Julia in the ``transverse'' direction. By contrast,  being the support of $\mu$, $J_2$ is meant to be ``repelling in all directions''.

In \cite{fs-hyperbolic} the authors show that is picture is indeed correct under (an adapted version of) the Axiom A assumption.

In the general case, one may at best expect that $J_1\setminus J_2$ is
filled with  Fatou disks in some measure-theoretic sense.   The  trace measure
$\sigma_T$ is a natural (non invariant) measure on $J_1$, and the question is whether a set of full trace measure in $J_1\setminus J_2$ is filled with Fatou disks. A related (stronger) problem is whether  $T\rest{J_1\setminus J_2}$ is a {\em laminar current} (as defined by Bedford, Lyubich and Smillie \cite{bls}).

This was shown to be true  in the basin of infinity for polynomial mappings of $\cd$ admitting an extension as a holomorphic mapping of $\pd$  by Bedford and Jonsson \cite{bj}. This question is also the main motivation in \cite{dt-genre} (see also \cite{dt-saddle}), where De Thélin gives some evidence for the  laminarity of $T$ outside $J_2$, and actually proves it for post-critically finite maps. The general case, however, remains open.

\medskip

In arbitrary dimension $k$, for $1\leq q\leq k$,  let $J_q = \supp(T^q)$. We then have a filtration  of the Julia set
\begin{equation}\label{eq:filtration}
 J=J_1\supset J_2\supset \cdots \supset J_k = \supp(\mu).
\end{equation}
Intuitively, the ``number  of  Fatou directions'' should decrease from $k-1$ on $J_1\setminus J_2$ to zero on $J_k$, so that one expects that  a set of full $\sigma_{T^q}$ measure of $J_q\setminus J_{q+1}$ is filled with Fatou disks of codimension $q$.

The list of well understood situations is even shorter in this case. For polynomial mappings of $\cc^k$ extending holomorphically to $\pk$,  it is shown in \cite{bj} that $T^{k-1}$ is laminar (with 1-dimensional leaves) in the basin of attraction of the hyperplane at infinity. General results about Fatou disks can be found in \cite{ueda, maegawa}.

\subsection{Fatou and Julia directions} Here we propose the following model for the various dynamical regimes along the Julia filtration \eqref{eq:filtration}, 
which is strongly reminiscent of the Oseledets multiplicative ergodic theorem.

\begin{conj}\label{conj:fatoumain}
Let $f$ be a holomorphic endomorphism of $\pk$ and $T$ be its Green current. Let $q$ be an integer with $1\leq q\leq k-1$. Then
for $\sigma_{T^q}$-a.e.  $x\in J_q\setminus J_{q+1}$
there exists a  complex sub-vectorspace $\mathcal{F}_x\subset T_x\pk$ of codimension $q$, such that:
\begin{enumerate}[i.]
\item if $v\in \mathcal{F}_x$, then $\displaystyle{\limsup_{n\cv\infty} \unsur{n}\log\norm{df^n_x(v)}\leq 0}$;
\item  if $v\notin \mathcal{F}_x$, then $\displaystyle{\limsup_{n\cv\infty} \unsur{n}\log\norm{df^n_x(v)}\geq  \frac{\log d}{2}}$.
\end{enumerate}
\end{conj}

Here $\norm{\cdot}$ refers to any norm on the tangent bundle.
A tangent vector $v\in T_x\pp^k$ such that
$\limsup_{n\cv\infty} \unsur{n}\log \norm{df^n_x(v)}\leq 0$ will be said
to be \textit{of Fatou type}. The set of such directions is a sub-vectorspace of
$T_x\pk$, called the {\em Fatou subspace}.
Thus, the conjecture asserts that $\sigma_{T^q}$-a.e. on $J_q\setminus J_{q+1}$, the Fatou subspace $\mathcal{F}_x$ has codimension $q$.

Notice that if $T^q$ were known to be laminar and filled with Fatou disks of codimension $q$ on $J_q\setminus J_{q+1}$, then to obtain {\em i.} it would be enough to consider the collection of tangent spaces to these Fatou disks.  Conversely, we do not  address the converse question whether the Fatou sub-bundle can be integrated to yield a laminar structure for $T^q$ on $J_q\setminus J_{q+1}$.

If we adopt the convention that $J_0 =\pk$ and $J_{k+1}=\emptyset$, then the conjecture is true for $q=0$ (by definition of the Fatou set) and $q=k$ (by the work of Briend-Duval \cite{briend-duval1}).
Thus item {\em ii.} may be seen as  a generalization in lower codimension
of the Briend-Duval bound on the Lyapunov exponents of $\mu$. 
Notice however that for $q\leq k-1$, $\sigma_{T^q}$ is \textit{not} an invariant measure, so {\em ii.} is not exactly a statement about Lyapunov exponents. Likewise, ergodic theoretic methods (like Oseledets-Pesin theory) do not seem well adapted to deal with this conjecture. 

One might also want to replace the limsup in {\em ii.} by a liminf.


\subsection{Results and methods} In this paper we prove several results towards this conjecture, including a complete proof for $q=1$. This, in particular,  settles the 2-dimensional case.

Our first main result is the following.

\begin{thm}\label{thm1}
Let $f$ be a holomorphic endomorphism of $\pk$ and $T$ be its Green current.
Let $q$ be an integer with $1\leq q\leq k-1$.

Then for $\sigma_{T^q}$-a.e.  $x\in J_q\setminus J_{q+1}$, the Fatou subspace $\mathcal{F}_x$ has dimension  {\em at least} $k-q$ at $x$.
\end{thm}

Simple examples show that the inequality in item {\em i.} of Conjecture \ref{conj:fatoumain} is not strict in general (see \S\ref{subsub:J12}). 

As already said, the trace measure $\sigma_{T^q}$ is not invariant so we cannot rely on ergodic theoretic methods  
 to prove this theorem. Instead, 
   we come back to the basic idea of a positive current  as being a   differential form with measure coefficients.
Let then  $S$  be an arbitrary positive closed current of bidegree $(q,q)$.
By duality, we can thus associate to such a $S$
  a measurable field of positive normalized $(p,p)$ vectors $t_S(x)$ ($p+q=k$),
and a positive measure $\sigma_S$ (the trace measure) such that the action of $S$ on test forms may be  expressed as
$$\bra{S,\phi}= \int \bra{t_S(x), \phi(x)} d\sigma_S(x).$$
This is known as the {\em integral representation of $S$}.
 If $t_S$ is well-defined at $x$,  we let $\mathrm{Span}(t_S(x))$ be the smallest subspace  $V\subset T_x\pk$ generating $t_S(x)$, i.e. such that
$t_S(x) \in \bigwedge_{(p,p)}V$. It is in a sense the tangent space to $S$ at $x$. Its
dimension can be any integer between $k-q$ and $k$, and will be referred to as the {\em rank} of $S$ at $x$. The directional information embedded in $S$ is   most  precise when its rank equals $k-q$. 
In Section \ref{sec:integral} we give several estimates on the rank of   general
positive (closed) currents. 

When $S$ is an invariant current, we obtain an invariant field of subspaces, hence a dynamically meaningful object. 
Albeit natural, it seems that this idea is used here for the first time in holomorphic dynamics. 

\medskip 

Back in the  context of endomorphisms on $\pk$, in Theorem \ref{thm:expansion} we give estimates on the expansion rate along the field of tangent spaces to the invariant currents $T^q$ (recall that $T$ is the Green current). In particular we obtain  that for generic $x\in J_q\setminus J_{q+1}$, $\mathrm{Span}(t_{T^q}(x))\subset \mathcal{F}_x$;  Theorem \ref{thm1} follows.
 An interesting point is that, in the event
 that $\sigma_{T^q}$ carries some mass on $\bigcup_{\ell>q}J_{\ell}$, we also obtain
estimates there. It turns out that the expansion rate along $T^q$ is never greater than $\frac{\log d}{2}$


\bigskip

In the second half of the paper, we study expansion properties in the
directions ``transverse to $T^q$'', motivated by assertion {\em ii.} of the conjecture.
For this, we extend to arbitrary dimensions a method
introduced by the author in \cite{birat} to obtain directional Briend-Duval type estimates for the Lyapunov exponents of birational mappings in dimension 2. This was adapted to  endomorphisms of $\pd$ by De Thélin \cite{dt-saddle}, and we follow his approach (see also \cite{ddg3} for related material).

Roughly speaking, the method is as follows. If $L$ is a generic $q$ dimensional linear subspace of $\pk$, consider the sequence of its iterates $f^n(L)$ and the associated  currents $S_n = d^{-qn} [f^n(L)]$.
Dinh \cite{dinh-lamin}   has shown that the major part
 of $f^n(L)$ has locally bounded geometry.
More precisely for every $\e>0$ there exists $r=r(\e)$ such that the part of
$f^n(L)$ which is {\em not} made of graphs of size $r$ over some direction
has mass less than $\e$ (relative to the trace measure of $S_n$).

If we  furthermore assume  that for generic $L$, 
most of the mass of $T^q\wedge S_n$ is concentrated on the bounded geometry part (hypothesis  $(H_q)$), then we show in Theorem \ref{thm:transversex} that Conjecture \ref{conj:fatoumain} holds in codimension $q$. For $q=1$, we are able to check this assumption (this essentially follows from \cite{isect}). This leads to the following (see Corollary \ref{cor:q=1}):

\begin{thm}\label{thm2}
Conjecture \ref{conj:fatoumain} is true for $q=1$.
\end{thm}

For general values of $q$ we   provide some strong evidence that 
the   assumption $(H_q)$ is always satisfied. We are also able to show that it holds (hence also Conjecture \ref{conj:fatoumain})  in the basin of a $q$-dimensional algebraic  (measure-theoretic) attractor (see Corollary \ref{cor:attractor}).

\subsection{Outline and acknowledgments} The plan of the paper is following.
In Section \ref{sec:positivity} we recall some preliminaries on positive exterior algebra. We then study
 in Section \ref{sec:integral}  the integral representation of positive closed currents in general. Endomorphisms of $\pk$  enter the picture in Section \ref{sec:tangentdyn}, where we study
  the rate of expansion of tangent vectors to the Green currents.

Section \ref{sec:geometry} is devoted to the study of the asymptotic geometry of varieties of the form $f^n(L)$, with $L$ a $q$-dimensional linear subspace. In particular we give some refinements of the  above mentioned    result of Dinh  and study the geometry of wedge products of the form $T^q\wedge [f^n(L)]$. In Section \ref{sec:transversex} we  turn this into  expansion results  transverse  to $T^q$.


\medskip

It is a pleasure to thank Henry De Thélin and Eric Bedford for helpful conversations.

\section{Preliminaries on positive exterior algebra}\label{sec:positivity}

\subsection{Vectors and covectors}
We start with some elementary considerations about vectors and covectors in the complex setting.

To avoid confusion, we  use the word covector for a form with constant coefficients, so that a diffential form is a field of covectors. Let $V$ be a complex vector space of dimension $k$. We use the notational convention that $p$ and $q$ are integers satisfying $p+q=k$.
By definition, $(p,p)$ vectors are the elements of the exterior algebra $\bigwedge^{p,p}V=\bigwedge^p V\otimes \bigwedge^p\overline V$, and
$(p,p)$ covectors are dual to $(p,p)$ vectors.
If $V$ is provided with
provided  with a basis $(e_i)_{i=1\ldots k}$,  then using  standard multi-index notation, a basis of the space of $(p,p)$ vectors is  $e_I\wedge \overline e_J$ with $\abs{I} = \abs{J}=p$.  By definition, $dz_I\wedge d \overline z_J$
is dual to $e_I\wedge \overline e_J$. Let us denote by ${}^\ast$ this duality. Finally, define $\langle \cdot, \cdot\rangle$ to be the associated $\cc$-bilinear  pairing, 
 which is normalized, by
declaring that $\langle{i^{p^2}e_I\wedge \overline e_J, i^{p^2}dz_I\wedge d \overline z_J}\rangle=1$ --this will ensure that pairing   positive objects results in nonnegative numbers.

The classical interpretation of the $p$ vector $u_1\wedge \cdots \wedge u_p$ in $\re^n$ is that of the sub-vector space  generated by $(u_1,\ldots, u_p)$ endowed with a volume form. In the complex setting, the geometric interpretation of $u\wedge \overline v$ is  a bit more cumbersome, because complex
conjugation here is understood as the involution sending (1,0) vectors to (0,1) vectors, not as the geometric conjugation of vectors in $V$, with respect to the almost complex structure. On the other hand,
$i u_1\wedge \overline u_1\wedge \cdots \wedge i u_p\wedge \overline u_p$ can indeed be interpreted as the complex sub-vector space  generated by $(u_1,\ldots, u_p)$ endowed with a volume form.

If $V$ is given a Hermitian metric, we let   $\beta$ be (twice) the associated (1,1) form. In coordinates,
 if $(e_1, \ldots , e_n)$ is any orthonormal basis, then $\beta = i\sum dz_j\wedge d\overline z_j$. 
With this convention   the volume form associated to the Hermitian metric is  $\frac{\beta^k}{2^kk!}$.
 We obtain an isomorphism $\Phi$ between $(p,p)$ vectors and $(q,q)$ covectors as follows:
if $t$ is a $(p,p)$ vector, we define $\Phi(t)$ as the unique $(q,q)$ covector s.t. for every $(p,p)$ covector
$\varphi$,
\begin{equation}\label{eq:dual}
\Phi(t)\wedge \varphi = \langle t,\varphi \rangle \frac{\beta^k}{k!}.
\end{equation} For instance,
$\Phi(i e_1\wedge \overline e_1\wedge \cdots \wedge i e_p\wedge \overline e_p ) = idz_{p+1}\wedge d\overline z_{p+1}\wedge \cdots \wedge idz_k\wedge d\overline z_k$.
Likewise, $\Phi$ may be expressed as  $\Phi(t) = \star t^\ast$ where $\star$ is the Hodge star.

\subsection{Positive (1,1) vectors and covectors}
We refer to Lelong  \cite{lelong} and Demailly \cite{demailly} for more details on the concept of positivity.
 Here we only gather some essential facts.
We work only with vectors, the case of covectors is completely similar.

Fix any basis $(e_1, \ldots, e_k)$ of $V$.  A $(1,1)$ vector
is positive (resp. real) if it writes as $t= i\sum_{i,j =1}^k t_{i,j} e_i\wedge \overline e_j$, with $(t_{i,j})$ a nonnegative Hermitian (resp.  Hermitian) matrix.
A positive $(1,1)$ vector is {\em decomposable} (or {\em simple}) if it may be written as $\lambda iu\wedge \overline u$, with $\lambda\geq 0$. Any positive $(1,1)$ vector is a sum of decomposable  ones. Let us denote by $P^{1,1}(V)$ (or $SP^{1,1}$ see below) the cone of positive (1,1) vectors.

The so-called {\em mass norm} on $(1,1)$ vectors is defined by $\norm{t} = \sum\abs{t_{i,j}}$.
Of course it depends on the choice of coordinates.

Assume  now that $V$ is equipped with a Hermitian metric and the basis $(e_i)$ is orthonormal.
 If $t$ is positive, the trace of the associated Hermitian matrix
 does not depend on the choice of an orthonormal basis, and is comparable to
 $\norm{t}$. To be specific, $\trace(t)\leq \norm{t} \leq k\trace(t)$.
 If $u\in \cc^k$ is of unit norm, then
$\mathrm{trace}(iu\wedge \overline u)=1$.

Define the \textit{rank} of a positive (1,1) vector $t$ to be the rank of the associated hermitian matrix. The rank can also be characterized as the least number of terms in a decomposition of $t$ as a sum of decomposable (1,1) vectors. In particular $t$ has rank 1 iff it is decomposable. Elementary considerations show that $\mathrm{rank}(t)$ is the largest integer $r$  s.t. $t^{\wedge r}$ is non-zero --from now on we write $t^r$ for $t^{\wedge r}$. In particular $t$ is decomposable iff $t^2=0$. We say that $t$ is strictly positive if $\rank(t)=k$ and degenerate if not. Then $t$ is degenerate iff $t^k=0$.

This discussion is valid  {\em mutatis mutandis} for (1,1) covectors. Recall from  \eqref{eq:dual} the duality
$\Phi$ between  $(k-1, k-1)$ vectors and  $(1,1)$ covectors. In particular we obtain that
a positive $(k-1, k-1)$ vector is decomposable iff $\Phi(t)^2=0$, and strictly positive iff
   $\Phi(t)^k\neq0$.

\subsection{Higher bidegree}\label{subs:higherdegree}
Recall that the space of $(k,k)$ vectors has dimension 1.
A $(k,k)$ vector is said to be positive if it is a nonnegative multiple of
$ie_1\wedge \overline e_1\wedge \cdots \wedge i e_k\wedge \overline e_k$.
A $(p, p)$ vector is {\em decomposable} if it is of the form
$ i u_1\wedge \overline u_1\wedge \cdots \wedge i u_p\wedge \overline u_p$.
By definition, a
{\em strongly positive} $(p, p)$ vector  is a convex combination of decomposable $(p,p)$ vectors. Let us  denote by $SP^{p,p}(V)$ the cone of strongly positive $(p,p)$ vectors on $V$.
A $(p,p)$ vector $t$ is {\em (weakly) positive}
 if for every strongly positive $(q,q)$ vector $t'$ ($p+q=k$),  the
$(k,k)$ vector $t\wedge t'$ is positive. It is also true that $t$ is strongly positive
if for any positive $t'$ of complementary degree, $t\wedge t'$ is positive. In other words, the cones of  positive and strongly positive vectors are dual to each other.  Following Lelong,
usually the single word ``positive"  is used as a shorthand for  ``weakly positive" (this also applies to currents).

 It is a fact that for $1<p<n-1$ the classes of positive  and strongly positive vectors differ. Notice that by definition strong positivity is  stable under wedge products (whereas weak positivity is not). In this paper we  mainly have to deal with strongly positive vectors\footnote{There is an alternate  notion of positivity which is intermediate between weak and strong, and gives rise to a self dual cone \cite{hk}.}.

\medskip

In standard multi-index notation, a  positive $(p,p)$ vector can be  written
as $t=i^{p^2} \sum t_{I,J} e_I\wedge \overline e_J$, where $t_{I,J} = \overline{t_{J, I}}$ and
$t_{I,I}\geq 0$ --for the $i^{p^2}$ see \cite[3.1.2]{demailly}. We define its {\em trace}   as $\mathrm{trace} (t)=\sum t_{I,I}$.  Notice that $\mathrm{trace} (t) = \langle t, \frac{\beta^p}{p!}\rangle$, so it does not depend on the choice of orthonormal coordinates.

It is clear that the dualities ${}^\ast$ and $\Phi$ send decomposable  positive vectors to decomposable
 positive covectors, hence they preserve strong positivity. Also, since ${t}^\ast\wedge {t'}^\ast = (t\wedge t')^\ast$, the duality ${}^\ast$ preserves weak positivity. Finally, from the property $\star\varphi\wedge \psi = \varphi\wedge \star\psi$ we conclude that $\Phi$ preserves positivity as well (see also \cite[pp. 64-65]{lelong}).

\medskip

We now discuss a notion of rank for positive $(p,p)$ vectors.
 If $t$ is a positive $(p,p)$ vector, we define $\mathrm{Span}(t)$ to be 
the smallest sub-vector space $W$ such that  $t\in \bigwedge^{p,p}(W)$, and $\rank(t)=\dim \mathrm{Span}(t)$.
We see that  $\rank(t)\geq p$ with equality iff $t$ is decomposable.

If $t$ is strongly positive, $t=\sum_{k=1}^s t_k$, where $t_k$ is (nonzero) decomposable,
then   $\mathrm{Span}(t)= \mathrm{Vect}(\mathrm{Span}(t_k),\ k=1\ldots s)$.

Also, still in case $t$ is strongly positive, $\rank(t)$ equals that of the (1,1) vector $t\llcorner\beta^{p-1}$. 
Indeed if $t$ is  decomposable,
$t=i e_1\wedge \overline e_1\wedge \cdots \wedge i e_p\wedge \overline e_p$, then
there exists $\lambda>0$ and an orthonormal family $(u_1, \ldots , u_p)$ such that
$\mathrm{Vect}(u_1, \ldots , u_p) =\mathrm{Vect}(e_1, \ldots , e_p)$ and
$t=\lambda i u_1\wedge \overline u_1\wedge \cdots \wedge i u_p\wedge \overline u_p$. We then infer that
$t\llcorner \beta^{p-1} = \lambda\sum_{j=1}^p iu_j\wedge \overline u_j$, hence the result.

%

 In view of applications to currents, we define the corank of a  $(q,q)$ covector $\phi$ to be the rank of the $(p,p)$ vector $\Phi^{-1}(\phi)$. Thus $\mathrm{corank}(\phi)\geq k-q=p$ with equality iff $\phi$ is decomposable.
An alternate characterization of the corank is given in \cite[p.65]{lelong}:  $k-\mathrm{corank}(\phi)$ is the greatest possible number of independent decomposable (1,1) forms dividing $\phi$, that is, independent linear forms $\alpha_1^\ast, \ldots, \alpha_r^\ast$ s.t. $\phi$ can be written as $\phi = i\alpha_1^\ast\wedge \overline\alpha_1^\ast\wedge \cdots \wedge i\alpha_{r}^\ast\wedge \overline\alpha_{r}^\ast\wedge \phi_1$.
%

\section{The integral representation of positive closed currents}\label{sec:integral}

\subsection{Preliminaries}
We just collect a few facts on positive currents
and again refer  the reader to \cite{lelong, demailly} for  details.
Since these notions are local we work in an open set $\om\subset \cc^k$.
As before, let $\beta ={i}\sum_{i=1}^k d z_i\wedge d\overline z_i$, and $p$, $q$ be integers with $p+q=k$.

A differential form of bidegree $(p,p)$ is (resp. strongly) positive if it satisfies this property at every point. A current $T$ of bidimension $(p,p)$ is (resp. strongly) positive if for every strongly positive (resp. weakly positive) $(p,p)$ test form $\varphi$, $\bra{T,\varphi}\geq 0$. Observe that  these notions are stable under weak convergence.

 It is well known that $T$ may be written in coordinates $(z_1, \ldots , z_k)$
as a $(q,q)$  form with measure coefficients
$T = {i^{q^2}} \sum T_{I,J} dz_I\wedge d\overline{z}_J$ where the $(T_{I,J})$ are complex  measures satisfying $T_{I,J} = \overline{T_{J,I}}$. Here the action of $T$ on a test form is expressed as $\varphi\mapsto \int T\wedge \varphi$

The {\em trace measure} of $T$ is defined by $\sigma_T = \sum_{\abs{I} = q} T_{I,I}$.
Notice  that $T\wedge \frac{\beta^p}{p!} = \sigma_T\frac{\beta^k}{k!}$.
Related to it is the {\em mass measure} $\norm{T} = \sum_{I,J}\abs{T_{I,J}}$. Throughout the paper,   notions of ``mass'' for positive currents will always be relative to $\sigma_T$.

There exists a constant depending only on $q$ such that
for all $I,J$, $\abs{T_{I,J}} \leq C_q\sigma_{T}$, so there exist measurable functions $f_{I,J}$ such that  $\abs{f_{I,J}}\leq  C_q$
and $T_{I,J} = f_{I,J}\sigma_T$.  Notice that $\sum f_{I,I} =1$. Thus we can
write
$T = \phi \sigma_T$, where
$\phi$ is a measurable field of positive $(q,q)$ covectors of trace 1.

Using the  duality \eqref{eq:dual}, we can formulate this by saying that
 there exists a measurable field $t_T$ of positive $(p,p)$ vectors of trace 1
 such that if $\varphi$ is any test $(p,p)$-form,
\begin{equation}
\label{eq:polar}
\langle T,\varphi \rangle = \int \langle t_T, \varphi \rangle \sigma_T
 \end{equation}
(from now on we omit the conventional $d$ of integration to avoid any confusion with exterior derivative or degree).
We  refer to either this representation or the representation $T=\phi\sigma_T$
 as the {\em integral representation} of $T$.

By the Lebesgue Density Theorem (see e.g. \cite[p.38]{mattila}), we can recover the tangent vector $t$ as follows: there exists a set $A\subset \supp(T)$ of full $\sigma_{T}$-mass such that if $x\in A$ and $\varphi$ is any test form, then
\begin{equation}\label{eq:lebesgue1}
\lim_{r\cv 0} \unsur{\sigma_{T}(B(x,r))}\langle T, \varphi\rest{B(x,r)}\rangle  =
\lim_{r\cv 0} \unsur{\sigma_{T}(B(x,r))} \int_{B(x,r)} \! \langle t(y), \varphi(y) \rangle \sigma_{T}(y) =
\langle t(x), \varphi(x)\rangle.
\end{equation}

Observe that  by continuity of $\varphi$, to ensure the existence of the limit in this equation,
 it is enough to test the convergence on constant forms $\varphi$.
Thus we can rewrite
 \eqref{eq:lebesgue1} as a convergence statement in the space of $(p,p)$ vectors (of trace 1)
 \begin{equation}\label{eq:lebesgue}
\lim_{r\cv 0} \unsur{\sigma_{T}(B(x,r))} \int_{B(x,r)} t(y)\ \sigma_{T}(y) =
 t(x).
\end{equation}

It would be interesting to investigate more precisely  the size
of the exceptional set $\supp(T)\setminus A$ for positive closed currents.
We do not address this problem here.


\begin{defi}
We say that a positive  current $T$ of bidimension $(p,p)$ is decomposable  at
 $x$ if the limit in \eqref{eq:lebesgue} exists  and $t(x)$ is a decomposable $(p,p)$ vector.
Likewise, if $t(x)$ is well defined, we define the rank of $T$ at $x$ to be $\rank(t(x))$.
If  $\rank(t(x))<k$
we say that $T$ is degenerate at $x$.
\end{defi}

Recall that $t(x)$ is   decomposable iff $\rank(t(x))=p$. For instance,
if $T$ is the integration current over a subvariety $M$ of dimension $p$, then it has rank $p$ a.e., and $t(x)$, which is well defined, at least at every smooth point, can be written as
$$t(x)=i\tau_1\wedge \overline  \tau_1\wedge \cdots \wedge i\tau_p\wedge \overline  \tau_p,$$ where $\tau_1, \ldots \tau_p$ is any orthonormal basis of $T_xM$ (the basis is not unique but $t(x)$ is).
As a consequence, laminar currents are decomposable a.e. \cite{bls}.

\subsection{Some notation}
It will be convenient  to us to work with the following
regularization procedure.
Let $\rho =\unsur{v_{2k}}\mathbf{1}_{B(0,1)}$ be the
 characteristic function of the unit ball normalized by its volume,
and consider the associated ``regularizing kernel'' $\rho_\e =
\unsur{\e^{2k}}\rho(\frac{\cdot}{\e}) = \unsur{v_{2k}\e^{2k}}\mathbf{1}_{B(0,\e)}$.
Given a positive current  $T$ we put $T_\e = T\ast \rho_\e =
\int \big((\tau_s)_*T \big) \rho_\e(s)ds$ which is positive and has continuous coefficients ($\tau_s$ is the translation of vector $s$). 
Likewise, if $\nu$ is a measure, we denote
$\nu\ast\rho_\e$ by $\nu_\e$. 
The mass of $\nu$ is denoted by $\m(\nu)$, and Lebesgue measure is denote by $\leb$.

Throughout the remaining part of this  section, if $T$ is a positive current of bidimension $(q,q)$, we let  $T= \phi \sigma_T$
(or $T= \phi_T \sigma_T$ when required) be its integral representation. Here $\phi$ denotes the associated field of $(q,q)$ covectors. We  only work with strongly positive currents, that is we require that $\phi_T$ is strongly positive a.e. Recall that this is not a restriction for $q=1$ and $q=k-1$.

\subsection{Pointwise self-intersections for absolutely continuous currents}

In this paragraph we consider a strongly positive current $T$, not necessarily closed, such that $\sigma_{T}$ is absolutely continuous w.r.t. Lebesgue measure. It has been observed by several
 authors, starting with \cite{bt} (see also \cite{boucksom}) that it is sometimes useful to work with the ``naive"  {\em pointwise self-intersection}
 $P(T^\ell)$, which is defined as follows. Write $T= \phi\sigma_T = \phi(x) h(x)dx$, and set
$P(T^\ell) = \phi^\ell h^\ell dx$. This is a differential form with Borel
coefficients, but not a priori a genuine current  since $h^\ell$ needn't be locally integrable.

On the other hand we have:
\begin{lem}
Assume that $T$ is strongly positive and $\sigma_{T}$ is absolutely continuous w.r.t. Lebesgue measure. If the family $(T_\e^\ell)$ has locally uniformly bounded mass as $\e\cv 0$ then
$P(T^\ell)$ is a well defined strongly positive  $(\ell q, \ell q)$ current.
\end{lem}

\begin{proof}
Write in coordinates $T= i^{q^2}\sum T_{I,J}dz_I\wedge d\overline z_J$ (resp.
$T_\e= i^{q^2}\sum (T_{I,J})_\e dz_I\wedge d\overline z_J$). The   Lebesgue Density Theorem implies
  that for all $I,J$, $(T_{I,J})_\e$ converges $\leb$-a.e. and in $L^1_{\rm loc}$ to $(T_{I,J})$.
Then by Fatou's Lemma we infer that
$$0\leq \int P(T^\ell) \wedge \beta^{k-\ell q} \leq \liminf\int  T_\e^\ell \wedge \beta^{k-\ell q} <+\infty,$$
so $P(T^\ell)$ has $L^1_{\rm loc}$ coefficients.
\end{proof}

In particular we see that if $(T^\ell_\e)$  converges in the sense of currents as $\e\cv 0$, its limit must be
$P(T^\ell)$. Therefore, if $T$ is a strongly positive current with absolutely continuous coefficients
and $T^\ell$ is well defined (in the sense that $T^\ell_\e$ converges to $T^\ell$)
it must coincide with $P(T^\ell)$. In this case  we may simply denote the pointwise
self intersection by $T^\ell$.

\medskip

Assume  now that  $T$ is an arbitrary strongly positive current.
It admits a Lebesgue decomposition
$T =  T_{\rm ac} + T_{\rm sing}$, induced by that of $\sigma_T$.
Notice that even when $T$ is closed, $T_{\rm ac}$ and $T_{\rm sing}$ are not closed in general. If
$(T_\e^\ell)$ has locally uniformly bounded mass as $\e\cv 0$
we can consider  $P(T^\ell_{\rm ac})$.
Again, when no confusion can arise we simply denote it by $T^\ell_{\rm ac}$. We conclude that
if $T^\ell$ is well defined, then $T^\ell_{\rm ac}\leq (T^\ell)_{\rm ac} \leq T^\ell$.

In the case $q=1$ and $\ell =k$, the pointwise self intersection has additional properties due to the concavity of
$M\mapsto (\mathrm{det}(M))^{1/k}$ in the cone of nonnegative Hermitian matrices \cite{bt, boucksom}.

\subsection{Wedge products} Here is our first main result on the integral representation of positive currents. It says  that the exterior powers of the tangent covectors can be read off the exterior powers of $T$. Notice that we do not assume $T$ to be  closed.

\begin{thm}\label{thm:wedge}
Let $T = \phi\sigma_T$
be a strongly positive current of bidegree $(q,q)$ in $\om\subset \cc^k$. Assume that for some
$\ell>1$, the family  $(T_\e^\ell)_{\e>0}$ has locally uniformly bounded mass as $\e\cv 0$.

Then $\phi^\ell(x) =0$ for $\sigma_T$-a.e. $x$ if and only if $T^\ell_{\rm ac}=0$.
\end{thm}

For $q=1$ and $\ell=k$ this result is somewhat implicit    in \cite[\S 5]{bt}.
The assumption on $T^\ell$ is always satisfied when $T$ is the restriction to $\om$ of a current on the  projective space $\mathbb{P}^k$. This is also true when $T^\ell$ is well defined in the sense of pluripotential theory. 

When $q=1$, since the decomposability of $\phi$ is detected by the vanishing of $\phi^2$ we immediately get the following corollary.

\begin{cor}
Assume that $T$ is a positive current of  bidegree (1,1) satisfying the assumptions of Theorem \ref{thm:wedge}. Then if  $\sigma_T\perp\leb$
 then $T$ is decomposable. Likewise, if $T^2$ is well-defined and $\sigma_{T^2}\perp\leb$
 then $T$ is decomposable.

If now the positive measure $T^k$ is well defined and $T^k\perp\leb$ then $T$ is degenerate a.e.
\end{cor}

Regarding higher values of $q$, the theorem only makes sense when $q\leq \frac{k}2$. Unfortunately, when $q>1$, the condition $\phi^2= 0$ does not impose severe restrictions on the rank of $T$ --and $\phi^3= 0$ does not restrict the rank at all.
 For instance if $\phi  = \sum_{j=2}^k idz_1\wedge d\overline z_1 \wedge idz_j\wedge d\overline z_j,$
we see that  $\phi^2 =0$ and the corank of $\phi$ is $k-1$.

\begin{cor}\label{cor:rank}
Assume that the hypotheses  of Theorem \ref{thm:wedge} are satisfied  with $\ell=2$ and arbitrary $q$. Then
if $T^2_{\rm ac}=0$,  $T$ has rank $<k$ a.e.
\end{cor}

\begin{proof}[Proof of the corollary]
We need to show that if $\phi$ is a strongly positive $(q,q)$ vector with $\phi^2=0$ then $\mathrm{corank}(\phi)<k$.
Write $\phi$ as a combination of decomposable $(q,q)$ covectors
$$\phi = \sum \lambda_\alpha i e_{1,\alpha}^\ast\wedge \overline  e_{1,\alpha}^\ast\wedge \cdots \wedge ie_{q,\alpha}^\ast\wedge \overline e_{q,\alpha}^\ast,$$ with $\lambda_\alpha>0$. Since $\phi^2=0$, for $\alpha$, $\alpha'$ in the decomposition we get that $\mathrm{Vect}(e_{1,\alpha}, \ldots e_{q,\alpha})\cap \mathrm{Vect}(e_{1,\alpha'}, \ldots e_{q,\alpha'})\neq \set{0}$. Thus there exists a nonzero
vector $u$ belonging to all these subspaces, and we infer that $iu^\ast\wedge \overline u^\ast$ divides $\phi$.
From the characterization of corank  given at the end of \S \ref{subs:higherdegree} we conclude that $\mathrm{corank}(\phi)<k$.
\end{proof}

\begin{proof}[Proof of the theorem]
We start with  an elementary lemma, whose proof is left to the reader. Recall that if $\nu$ is a measure in $\re^d$,  $\nu_\e$ stands for $\rho_\e\ast\nu$.
\begin{lem}
Let $\nu$ be a measure in $\re^d$, and $f\in L^1(\nu)$.
Then $(f\nu)_\e = f_\e\nu_\e$, where $$f_\e(x) = \unsur{\nu(B(x,\e))}\int_{B(x,\e)} f\nu,$$ and 
$\nu_\e$ is absolutely continuous, with
$$\frac{d\nu_\e}{dx}= \frac{\nu(B(x,\e))}{v_{d}\e^{d}}.$$
\end{lem}
%


From  this lemma we can relate the integral representations of $T$ and $T_\e$. Write
$T= \phi\sigma_T$. Notice first that since $\sigma_T = T\llcorner \frac{\beta^p}{p!}$ and $\beta$ is translation invariant, we have that $\sigma_{T_\e} = (\sigma_T)_\e$. Then, using the  above lemma,
and denoting $T_\e = \phi_\e(\sigma_T)_\e$, we get that
$$\phi_\e(x) = \unsur{\sigma_T(B(x,\e))} \int_{B(x,\e)}\phi(y)\sigma_T(y).$$
By  Lebesgue Density \eqref{eq:lebesgue}, $\phi_\e(x)$ converges to $\phi(x)$ a.e.


Now assume that there is a set $E$ of positive trace measure such that $\phi^\ell$ is nonzero on $E$. To prove the theorem we will show that $\sigma_T\rest{E}$ is absolutely continuous, thus $T\rest{E}\leq T_{\rm ac}$. It then immediately follows that $T^\ell_{\rm ac}\rest{E}$ is nonzero.
Observe that by definition of $T^\ell_{\rm ac}$, the converse implication is obvious.
Let $E_c$ be  the set of $x$ s.t. $\phi^\ell(x)\wedge{\beta^{k-\ell q}}\geq c{\beta^k}$. For small enough $c>0$, $\sigma_T(E_c)>0$.
It is enough to prove that $\sigma_T\rest{E_c}$ is absolutely continuous.

Let $R = T\rest{E_c} = \phi_R\sigma_R$. The family $(R_\e)$ has locally uniformly bounded mass as $\e\cv 0$.
Write $R_\e =\phi_{R_\e}\sigma_{R_\e}$; since $\sigma_{R_\e}$ is absolutely continuous
we  define $h_\e \in L^1_{\rm loc}$ by $\sigma_{R_\e} = h_\e  \leb$. The proof will be finished if we
 show that $(h_\e)_{\e>0}$ is locally bounded in $L^\ell$. Indeed let then $h$ be a cluster value of this family relative to the weak* topology in $L^\ell$. Since $\sigma_{R_\e}\cv\sigma_R$ weakly as measures, we infer that $\sigma_R= h\leb$, hence the result.

Let us prove our claim. We have
that $R_\e^\ell = \phi_{R_\e}^\ell h_\e^\ell \leb$, with $$\phi_{R_\e} (x)
=\unsur{\sigma_R(B(x,\e))} \int_{B(x,\e)}\phi_R(y)\sigma_R(y), $$ and for every
$y$, $\phi_R^\ell(y)\wedge {\beta^{k-\ell q}}\geq c{\beta^k}$. By Lemma
\ref{lem:caratheodory} below there is a $\delta>0$ such that $\phi_{R_\e}^\ell\wedge \beta^{k-\ell q}\geq c\delta \beta^k$, from which we infer that
$$ O(1) = \int R_\e^\ell \wedge \beta^{k-\ell q} \geq c\delta \int
  h_\e^\ell {\beta^k},$$ which was the desired estimate.
\end{proof}

\begin{lem}\label{lem:caratheodory}
Let $(\phi_\alpha)_{\alpha\in \mathcal{A}}$ be a measurable family of strongly
positive $(q,q)$ covectors of trace 1, s.t. for each $\alpha$, $\phi_\alpha^\ell \wedge {\beta^{k-
\ell q}}\geq c{\beta^k}$. Let $\nu$ be a probability measure on $\mathcal{A}$
and $\phi  = \int \phi_\alpha d\nu(\alpha)$. Then there is a constant $\delta$
depending only on $k$, $q$, and $\ell$  such that $\phi^\ell \wedge
{\beta^{k-\ell q}}\geq  c\delta \beta^k$.
\end{lem}

\begin{proof}
$\phi$ belongs to the closed convex hull of $(\phi_\alpha)_{\alpha\in \mathcal{A}}$.
Without  loss of generality, we may assume that the family $(\phi_\alpha)$ is closed, and it is also bounded because the set of strongly positive $(q,q)$ covectors of trace 1 is.
Thus we conclude that $\phi$ belongs to the  convex hull of $(\phi_\alpha)_{\alpha\in \mathcal{A}}$. By Caratheodory's Theorem there is a finite subset $\set{\phi_{\alpha_i}, \ i= 1\ldots d+1}$ (where $d=\left(\begin{smallmatrix} q\\k\end{smallmatrix}\right)^2-1$ is the dimension of the ambient affine space) such that
$\phi$ belongs to the convex hull of the $\phi_{\alpha_i}$. We conclude that
$$\phi^\ell\wedge{\beta^{k-\ell q}}= \left( \sum \lambda_i\phi_{\alpha_i}\right)^\ell\wedge{\beta^{k-\ell q}} \geq \sum \lambda_i^\ell\phi_{\alpha_i}^\ell\wedge{\beta^{k-\ell q}}\geq c \left( \sum \lambda_i^\ell\right){\beta^{k}}\geq c (d+1)^{1-\ell} {\beta^{k}},$$ where the first inequality follows from the fact that a product of strongly positive covectors is strongly positive, and the last one from H\"older's inequality and the fact that $\sum\lambda_i=1$.
\end{proof}

From the dynamical point of view, here is an interesting open question: what sort of relationship is there between the integral representations of $T$ and $T^q$, when  $T$ is a positive closed current of bidegree $(1,1)$?
A basic difficulty here  is that in general $\sigma_{T^q}$ and $\sigma_T$ are mutually singular.

\subsection{Projections}

Here we show  that on projective space,
Theorem \ref{thm:wedge} together with a projection argument leads an interesting estimate on the rank of positive currents of bidegree $(q,q)$ with $q>1$. Since we do not use it in the sequel, we do not include the proof (see \cite{noterank} for details).
  
Fix a Fubini study metric and denote the associated K\"ahler form by $\omega$. From now on the notions of trace, etc. will be relative to this metric.

The {\em dimension}  of a measure $\mu$ is defined as
$$\dim(\mu) =\inf \set{\mathrm{HD}(E), \ E \text{ Borel set with } \mu(E)=1},$$
where $\mathrm{HD}$ denotes Hausdorff dimension.
 If $T$ is a positive closed current of bidimension $(p,p)$, then
$\dim (\sigma_T)\geq 2p$.

\begin{thm}\label{thm:proj}
Let $T$ be a strongly positive closed current of bidimension $(p,p)$ on $\mathbb{P}^k$, and assume that
$\dim(\sigma_T)< 4p$.
Then $\sigma_T$ a.e. we have that
\begin{equation}\label{eq:dim}
p\leq \mathrm{rank} (T) \leq \unsur{2}\dim  (\sigma_T).
\end{equation}
\end{thm}

For instance, if $\dim  (\sigma_T)< 2(p+1)$, then $T$ is decomposable a.e.
The assumption on $\sigma_T$ is void if $p> k/2$. In general it is unclear whether it is necessary.
The estimate \eqref{eq:dim} 
 is sharp since if $r\geq 0$ and
$V$ is any linear subspace of dimension $p+r$, $T= [V]\wedge \omega^r$ is a positive closed current of bidimension $(p,p)$ with rank $p+r$  everywhere and $\dim  (\sigma_T)=2(p+r)$.

%

\section{Upper estimates for tangential expansion}\label{sec:tangentdyn}

\subsection{Preliminary considerations} We start with some classical facts. We refer the reader to the survey papers \cite{sibony, guedj-survey, ds-survey} for more details and references.
Let $f$ be a holomorphic endomorphism of $\pp^k$ of degree $d>1$. Recall that $f$ is given in homogeneous coordinates by $k+1$ homogeneous polynomials of degree $d$ without non-trivial common zero. Fix a Fubini-Study metric $\norm{\cdot}$ on $T\pk$ with associated $(1,1)$ form $\omega$. If $X\subset\pk$
is an analytic subset, we use the notation $\omega_X$ for the restriction of $\omega$ to $X$.
The {\em Green current} of $f$ is defined as $T = \lim_{n\cv\infty} \unsur{d^n} (f^n)^*\omega$.

More precisely, let $g_1$ be the continuous quasi-psh function  defined by $\unsur{d}f^*\omega = \omega + dd^cg_1$ and $\max g_1 = 0$. Then
$\unsur{d^n} (f^n)^*\omega = \omega + dd^c g_n$, where $g_n = \sum_{k=0}^{n-1} \unsur{d^k} g_1\circ f^k$,
hence $g_n$ converges uniformly to the  quasi-psh function $g= \sum_{k=0}^\infty
\unsur{d^k} g_1\circ f^k$, and  $T = \omega+ dd^cg$.

Likewise, if $H$ is a generic hyperplane,  $\lim_{n\cv\infty} \unsur{d^n} (f^n)^*[H] = T$.

\medskip

Recall  that the {\em $q$-th Julia set} $J_q$ is defined by $J_q= \supp(T^q)$. $J_1$ is the Julia set of $f$ in the ordinary sense. It follows from intersection theory of currents that $T^q =\lim_{n\cv\infty} \unsur{d^{nq}} (f^n)^*(\omega^q)$.
Consequently, being a weak limit of strongly positive forms, $T^q$ is strongly positive.

\begin{defi}
A tangent vector $v\in T_x\pp^k$ is said
to be \textit{of Fatou type} if $$\limsup_{n\cv\infty} \unsur{n}\log \norm{df^n_x(v)}\leq 0.$$
The Fatou subspace $\mathcal{F}_x\subset T_x\pp^k$ is the sub-vector space of Fatou directions.
\end{defi}

The collection of  Fatou subspaces defines
 a forward and backward invariant subbundle: $df_x(\F_x)\subset \F_{f(x)}$ and  $df_{f(x)}^{-1}(\F_{f(x)})\subset \F_{x}$.

The term `Fatou' is convenient but somewhat misleading, since of course the definition  does not prevents from subexponential expansion. Notice that even if we replace it by the stronger condition that $\norm{df^n_x(v)}$ is bounded, every tangent vector at an indifferent periodic point is Fatou, while indifferent periodic points themselves may belong to either the Fatou or the Julia set.
On the other hand the relevance of this definition is partly justified by the following dichotomy, which suggests that Julia-like behavior is always related to some exponential growth of the derivative.

\begin{prop}
Assume that $\Delta$ is a one-dimensional holomorphic disk in $\pp^k$. Then
\begin{itemize}
\itm either $(f^n\rest{\Delta})_{n\geq 0}$ is a normal family;
\itm or $\displaystyle{\liminf_{n\cv\infty} \unsur{d^n}\vol(f^n(\Delta))>0}$.
\end{itemize}
\end{prop}

\begin{proof}
Forn\ae ss and Sibony showed \cite[Prop. 5.10]{fs2} that
 that  $T\wedge [\Delta] = 0$ iff $(f^n\rest{\Delta})_{n\geq 0}$ is a normal family.
So it is enough to prove that if $\m(T\wedge [\Delta])>0$, then $\liminf \unsur{d^n}\vol(f^n(\Delta))>0$, which is known to be true (see e.g. \cite[Prop. 5.3]{fs-hyperbolic}; the result is stated in dimension 2 there but the adaptation to the general case is obvious).
\end{proof}

A theorem of Berteloot and Dupont \cite{berteloot-dupont} asserts that $\mu\ll \leb$ if and only
if $f$ is a {\em Lattès example}, that is,  a quotient of  a linear map on a complex torus. Hence,
from Theorem \ref{thm:wedge} we  immediately get the following corollary.

\begin{thm}\label{thm:cor}
Let $f$ be an endomorphism of $\pp^k$ of degree $d>1$. Then the Green current is degenerate a.e. (decomposable a.e. if $k=2$) unless $f$ is a Lattès example.

Furthermore, $T$ is always decomposable a.e.  on $J_1\setminus J_2$
\end{thm}

This result means that, except in the case of Lattès examples, the Green current always carries some directional information.
Theorem \ref{thm:wedge} also admits  consequences on the tangent vectors of  $T^q$ for $1<q\leq \frac{k}{2}$; we leave the precise formulation to the sagacity of the reader. 

\subsection{Expansion for tangent vectors}
On the $q^{\rm th}$ Julia set $J_q$, the invariant current $T^q$ induces a measurable $df_*$-invariant subbundle $\mathcal{T}^q$ of $T\pk$,
whose stalk is defined at $\sigma_{T^q}$ a.e. $x$  by  $\mathcal{T}^q_x =
\mathrm{Span}(t_{T^q}(x))$.
In particular $\dim \mathcal{T}^q_x\geq k-q$, with equality iff $T^q$ is decomposable at $x$. 
The invariance of $\mathcal{T}^q$ simply follows from the fact that $f$ is a local diffeomorphism a.e., since $T^q$ gives no mass to the critical set.

The following result describes the
expansion properties of the  action of $f$ on this invariant subbundle. Its sharpness will be discussed in \S \ref{subs:comments}.

\begin{thm}\label{thm:expansion}
Let $f$ be a holomorphic endomorphism of $\mathbb{P}^k$, and $T$ be its Green current.
\begin{enumerate}
\item[i.] For  $\sigma_{T^q}$-a.e. $x\in J_q$, if $v\in \mathcal{T}^q_x$, then $$\limsup_{n\cv\infty} \unsur{n}\log \norm{df^n(v)}\leq \frac{\log d}{2}.$$
\item[ii.]For  $\sigma_{T^q}$-a.e. $x\in J_q\setminus J_{q+1}$, if $v\in \mathcal{T}^q_x$, then $$\limsup_{n\cv\infty} \unsur{n}\log \norm{df^n(v)}\leq 0.$$
\end{enumerate}
\end{thm}

\begin{cor}\label{cor:fatouminor}
For  $\sigma_{T^q}$ a.e. $x\in J_q\setminus J_{q+1}$, $\mathcal{T}^q_x\subset \mathcal{F}_x$. In particular the Fatou subspace has dimension $\geq k-q$ at $x$.
\end{cor}

Thus we have produced at least $k-q$ Fatou directions a.e. on $J_q$.
On the other hand nothing  prevents {\em a priori} $\mathcal{F}_x$ from being larger than $\mathcal{T}^q_x$.
For instance this happens at indifferent periodic points.  Upper estimates on  the dimension of $\mathcal{F}_x$  will be given in Section \ref{sec:transversex}.

\begin{proof}
Let $E$ be the set of points   $x\in J_q$ such that the tangent vector $t_{T^q}(x)$ is well defined.
Then $\sigma_{T^q}(E)=1$ and if $x\in E$,  $t_{T^q}(x)$ is a  strongly positive $(k-q, k-q)$ vector of rank $\geq k-q$. Consider the $(1,1)$ vector $t_{T^q}(x)\llcorner\omega(x)^{k-q-1}$: it  of the form $\sum_j \lambda_j iu_j\wedge \overline u_j$, where $(u_j)$ is an orthonormal basis of $\mathcal{T}^q_x=\mathrm{Span}(t_{T^q}(x))$ and $\lambda_j>0$.
The numbers $\lambda_j$ are intrinsic since they are the eigenvalues of the Hermitian matrix associated
to $t_{T^q}(x)\llcorner\omega(x)^{k-q-1}$. Since $\trace (t_{T^q}(x))=\bra{t_{T^q}(x), \omega^{k-q}}=1$, $\sum\lambda_j=1$. On the other hand, the $\lambda_j$ needn't be bounded below.  For $c>0$ define  $E_c = \set{x\in E,\
\min\lambda_j\geq c}$; we have that $\bigcup_{c>0}E_c = E$, so it is enough to prove item
{\em i.} for $x\in E_c$.

Let now $x\in E_c$ and let $v\in \mathcal{T}^q_x$ be of unit norm. We have
\begin{align*}
\norm{df_x^n(v)}^2 &\leq \trace \lrpar{(df^n_x)_* (iv\wedge \overline v)}\leq \unsur{c}
\trace  \lrpar{(df^n_x)_* (t_{T^q}(x)\llcorner\omega(x)^{k-q-1})} \\
 & =\unsur{c}\bra {(df^n_x)_* (t_{T^q}(x)\llcorner\omega(x)^{k-q-1}), \omega(f^{n}(x))}\\
 &= \unsur{c}\bra{t_{T^q}(x),(df^n_x)^*(\omega(f^{n}(x)))\wedge\omega(x)^{k-q-1}}
\end{align*}
(recall that $df^n$ is invertible a.e.). Integrating this expression we obtain that
\begin{align*}
\int_{E_c} \norm{df_x^n\rest{\mathcal{T}^q_x}}^2\sigma_{T^q}(x)&\leq
\unsur{c}\int\bra{t_{T^q}(x),(df^n_x)^*(\omega(f^{n}(x)))\wedge\omega(x)^{k-q-1}} \sigma_{T^q}(x)\\
&= \unsur{c}\int T^q\wedge (f^n)^*\omega\wedge \omega^{k-q-1} = \unsur{c}d^n,
\end{align*}
where the last equality comes from cohomology. We can now finish the proof of {\em i}  by using the
Borel-Cantelli lemma:
let $E_{c, n} = \set{x\in E_c, \ \norm{df_x^n\rest{\mathcal{T}^q_x}}^2 \geq n^2d^n}$,
we infer that $\sigma_{T^q}(E_{c, n}) \leq \frac{c}{n^2}$. Thus a.e. $x$ belongs to finitely many
$E_{c, n}$'s and we are done.

\medskip

The proof of {\em ii.} is similar. Define now $F_c=\set{x\in E_c,  \ \dist(x, J_{q+1})\geq c}$. Introducing
$F_{c, n} = \set{x\in F_c, \ \norm{df_x^n\rest{\mathcal{T}^q_x}}^2 \geq n^2}$
 and using the Borel-Cantelli lemma again, the proof will be finished if we show that
$\int_{F_c}T^q\wedge (f^n)^*\omega\wedge \omega^{k-q-1}$ is uniformly bounded in $n$.
For this, let $\chi$ be a cut-off function, with $\chi=1$ on $F_c$ and  $\chi=0$  on $J_{q+1}$,
and write 
\begin{align*}
\int_{F_c}T^q\wedge (f^n)^*\omega\wedge \omega^{k-q-1} &\leq
\int \chi T^q\wedge (f^n)^*\omega\wedge \omega^{k-q-1}\\
&=\int \chi T^q\wedge ((f^n)^*\omega-T) \wedge \omega^{k-q-1} \text{ (because }T^{q+1}=0 \text{ on }\supp(\chi))\\
&= d^n \int \chi T^q\wedge dd^c(g_n-g) \wedge \omega^{k-q-1} \\
&= d^n \int (g_n-g) T^q\wedge dd^c\chi\wedge \omega^{k-q-1} \\
&= O(1)\text{ (since by construction }\abs{g_n-g}=O(d^{-n})),
\end{align*}
which was the desired estimate.
\end{proof}

\subsection{Examples and comments}\label{subs:comments}
In this paragraph we illustrate Theorem \ref{thm:expansion} with several examples. We 
restrict to the   2-dimensional case, which  is already quite rich.

\subsubsection{} A first possibility is that $f$ is a Lattès example. Then $T$ is strictly positive a.e. and $\sigma_T$ as well as $T\wedge T=\mu$ are absolutely continuous w.r.t. Lebesgue measure. In this case the inequality in Theorem \ref{thm:expansion} {\em i.} is an equality, since by the work of Briend and Duval \cite{briend-duval1} we know that the Lyapunov exponents of $\mu$ are never smaller than $\frac{\log d}{2}$. This actually yields a new proof of the minimality  of the Lyapunov exponents of Lattès examples.

\medskip

If $f$ is not Lattès, $T$ is decomposable a.e., so it contains directional information.
An interesting situation is when $\mu \ll\sigma_T$ but $\mu\perp\leb$.
This phenomenon happens for instance for mappings of the form
$[P(z,w):Q(z,w):t^d]$ on $\pd$, with $[P:Q]$ a Lattès example on $\pu$. Indeed, working on $\set{t\neq 0} \simeq \cd$, if we let $B$ be the basin of attraction of 0, $\fr B$ is locally spherical outside a set of the form $\pi^{-1}(C)$, where $\pi:\cd\setminus 0\cv\pu$ is the natural map and $C$ is a finite set (see \cite{berteloot-loeb1} for this and \cite{dupont-spherique} for similar results in higher dimension).
On this locally spherical part, $\mu$ and  $\sigma_T$ are absolutely continuous with respect to the natural area measure --the local structure of $T$ is that of $\log^+\norm{z}$.

\medskip

In this case we have the following.

\begin{thm}
 Let $f$ be a holomorphic endomorphism of $\pd$ of degree $d\geq 2$, 
 such that  $\mu \ll\sigma_T$ but $\mu\perp\leb$, where $T$ is the Green current and $\mu = T^2$ is the equilibrium measure.

Then $\mu$ has Lyapunov exponents $\displaystyle \chi_2 > \chi_1 = \frac{\log d}{2}$.
\end{thm}

\begin{proof} 
By \cite{briend-duval1}, the Lyapunov exponents satisfy $\chi_2 \geq  \chi_1 \geq \frac{\log d}{2}$.
 We know from the work of Berteloot and Dupont \cite{berteloot-dupont}  that $\chi_2$ must be greater than $\frac{\log d}{2}$ for otherwise $f$ would be a Lattès example and $\mu$ would be absolutely continuous w.r.t. Lebesgue measure. 
On the other hand, Theorem \ref{thm:expansion} provides an invariant field of directions along which the expansion rate is not greater than $\frac{\log d}{2}$. The result follows.
\end{proof}

This raises the following interesting question.

\begin{question}\label{q1}
 Is the converse true? That is, if the minimal exponent of $\mu$ is $\frac{\log d}{2}$, does one have $\mu\ll\sigma_T$?

Is this  a rigid situation? That is, if  $\mu\ll\sigma_T$, does $f$ admit a 1-dimensional Lattès factor in some sense? Is $f$ a Lattès-like mapping in the sense of \cite{favre-pereira}?
\end{question}

\medskip

When $\mu$ and $\sigma_T$ are mutually singular, it is still possible that $\sigma_T(J_2)>0$. An extreme instance of this happens when $J_2=\pd$. Then we have an invariant field of complex lines  defined $\sigma_T$-a.e. on $J_2$ which may be expanded by the dynamics.

To get a simple example, consider a non-Lattès rational map $h$ on $\pu$  possessing  an ergodic
 measure $\nu$ of positive Lyapunov exponent equivalent to Lebesgue measure
 (this exponent is smaller than $\frac{\log d}{2}$ by the Pesin formula). This phenomenon occurs on a set of positive measure on the space of rational maps, including all critically finite maps without superattracting cycles \cite{rees} (see also  \cite{eremenko-lyubich} for an account on this). Let $\widehat f: \pu\times \pu \circlearrowright$ be defined by $\widehat f(z,w) = (h(z), h(w))$.  Taking the quotient of $\pu\times \pu$ by $(z,w)\sim(w,z)$,
 $\widehat f$ descends to a holomorphic map on $\pd$, whose properties can easily be read-off from
 $\widehat f$ \cite{sibony}. In particular the Green current viewed on on $\pu\times\pu$ is $\pi_1^*\mu_h + \pi_2^*\mu_h$, where $\pi_j:\pu\times\pu\cv\pu$ are the natural projections
Geometrically speaking, $\widehat T$ is a uniformly woven current (see below \S\ref{subs:prelwoven}) of the form
$$\int[\pu\times\set{w}] d\mu_h(w) + \int[\set{z}\times \pu] d\mu_h(z), $$ whose trace measure is $\leb\otimes\mu_h + \mu_h\otimes \leb$. From this we easily obtain that the generic expansion rate along $T$  is given by the Lyapunov exponent of $\nu$.

\subsubsection{}\label{subsub:J12} Let us now discuss the dynamics on $J_1\setminus J_2$. The field of tangent vectors to $T$ induces  a measurable 1-dimensional  invariant sub-bundle of $T\pd$, contained in the Fatou sub-bundle $\mathcal{F}$ -we'll see in Corollary \ref{cor:q=1}  that the two actually coincide. Since $\sigma_T$ is not invariant, the sequence $\norm{(df^n)_*t_T(x)}$ needn't  converge. Nevertheless we can define two measurable  invariant subsets
$$E^-= \set{x, \limsup\unsur{n} \log \norm{(df^n)_*t_T(x)}<0}\text{ and }E^0 = 
\set{x, \limsup\unsur{n} \log \norm{(df^n)_*t_T(x)}=0}.$$ 
In this respect, quotients of mappings $(h(z),h(w))$ on $\pu\times\pu$  
are easy to analyze. A Siegel disk for $h$ (resp. an attracting basin) gives rise to  a region where the expansion rate is zero (resp. negative) --this example appears in \cite{dt-saddle}. We see in particular that both $E^-$ and $E^0$ can simultaneously be of   positive trace measure.

In a slightly different fashion, in \cite{dt-saddle} De Thélin studies invariant measures of the form
$\nu=T\wedge S$, where $S$ is a cluster value of push-forwards of lines: 
$S = \lim_{j\cv\infty} \unsur{n_j}\sum_{k=1}^{n_j} \unsur{d^k}f^k_*[L]$. He shows that $\nu$ admits a positive Lyapunov exponent, and that if furthermore $\nu(J_1\setminus J_2)>0$ 		and $\nu$ carries no mass on analytic subsets, it admits a nonpositive exponent, that is, it is of (weak) saddle type.

It is natural to try to relate these saddle measures and contraction properties along $T$. Here is a specific question:

\begin{question}\label{q2}
Assume that such a saddle measure admits a negative exponent. Is $E^-$ of positive trace measure? Is $T$ laminar (i.e. described by stable manifolds) there?

Conversely, 
does contraction on a set of positive trace measure implies the existence of a 
saddle measure with negative exponent? 
\end{question}

For mappings that are polynomial on $\cd$, the work of Bedford-Jonsson \cite{bj} provides a satisfactory answer to the first part of the question (here the saddle measure is supported on the line at infinity). The situation is also well understood for mappings satisfying Axiom A \cite{fs-hyperbolic}. Besides these cases, the problem is open, even in the presence of an attractor--see e.g. the questions in \cite[\S 6]{dinh-attractors}.

%

\section{Geometry of iterated subvarieties}\label{sec:geometry}

From now on, we study the expansion properties of the dynamics in the directions ``transverse to'' the current $T^q$. In the next section we will give some sufficient conditions for expansion depending on the geometry of certain iterated subvarieties. Here we formalize the idea of a subvariety having bounded geometry on a certain subset and  give explicit bounds for the geometry of $f^n(L)$, where as above $f$ is an endomorphism of $\pk$ and $L$ is a linear subspace.

\subsection{Definitions and preliminaries}\label{subs:prelwoven}

\subsubsection{}\label{subsub:bddgeom} Let $V$ be an analytic subset of $\pk$, of pure dimension $q$. We say that $V$ has {\em bounded geometry at scale $r$} at $x$ if $V$ contains  a graph through $x$, of diameter bounded by $D$, over the ball of radius $r$ in its tangent space at $x$ (relative to the orthogonal projection on $T_xV$). We denote by $V[r]$ the set of such points $x$.  Here $D$ is a constant that may be fixed freely (say,  small as compared to the diameter of a coordinate patch) It is also understood that $r$ is small with respect to $D$.  In the situation where $V$ is locally reducible at $x$ we take the union of possible graphs.  In particular   $V[r]$ may have (mild) singularities.

Our aim here is to estimate how close $V[r]$ is to $V$ when $r$ is small, for certain classes of dynamically defined varieties $V$. 

This idea is intimately related to the theory of geometric (that is woven and laminar) currents, that is, currents that are integrals of varieties outside a set of arbitrary small mass. Indeed, if $V_n$ is a sequence of varieties, of volume $v_n$, such that for fixed $r$, $\vol(V_{n}[r])\geq v_n(1-\e(r))$,
with $\e(r)\cv 0$ as $r\cv 0$ {\em uniformly in $n$}, then the cluster values of the sequence of currents
$v_n^{-1}[V_n]$ have some geometric structure, since  the bounded geometry part passes to the limit. This theory was developed in e.g. \cite{bls, lamin, dt-boule, dinh-lamin, dt-lamin}, and has many dynamical applications. The situation in the present paper is slightly different, since we only need to estimate the geometry for large, but finite, $n$, and do not  need to consider the limiting objects. In particular $r$ may be allowed to decrease to zero as $n$ tends to infinity.

\subsubsection{}\label{subsub:cube} A technically convenient way to understand $V[r]$ is the following. Assume that we are working in a ball inside some coordinate chart, and all subdivisions, etc., are relative to this ball. 
Notice first that by the Cauchy estimates, there exists a constant $K$ such that 
if $V$ contains a graph through $x$ over $B(\pi(x), Kr)$, relative to {\em some} orthogonal projection $\pi$, then $x\in V[r]$. Thus we can fix once for all a family of
 $\lrpar{\begin{smallmatrix}
q\\k \end{smallmatrix}}$ orthogonal projections $\pi_j$ to $\cc^q$ in general position, and look for graphs over these projections. 

Given one of the $\pi_j$, 
we consider a subdivision $\qq_j$ of the projection base by cubes of size $r$. If $Q\in \qq$, we declare that an {\em irreducible} component  of $\pi_j^{-1}(Q)\cap V$ is {\em good} if it is a graph over $Q$ of diameter $\leq D$, {\em bad} otherwise. We usually denote by $V_\qq$ the union of good components.

We also obtain a subdivision $\mathcal{C}$ of $\cc^k$ by affine cubes of size $O(r)$, whose atoms are the $\bigcap \pi^{-1}_j(Q_j)$, $Q_j\in\qq_j$.
By taking the union of the good components relative to all projections, we get a variety $V_{\mathcal{C}}$, which in each cube is a union of graphs relative to the $\pi_j$. There exists a constant $K$ depending only on the  $\pi_j$ such that $V[Kr]\subset V_\cC$ (apart from points at the boundary of the cubes of $\cC$).

Conversely,  if $Q$ is a cube of size $r$, and
$\lambda>0$, we let $\lambda Q$ be the cube with the same center as $Q$, 
  homothetic to it by a  factor $\lambda$. 
If $\pi$ is a projection as above,  a  {\em strong good component} over $Q$ 
is by definition the restriction to of $\pi^{-1}(Q)\cap V$  of a good component over 
$2Q$. We can   construct $V_{\qq, \rm strong}$, and $V_{\cC, \rm strong}$ in the same way as before, by keeping only strong good components, and  infer that for some $K'$, $V_{\cC, \rm strong}\subset V[K'r]$.

Another way to present this construction is the following: consider  for each projection a family of overlapping subdivisions by cubes of size $2r$, with the property that every point in the base  is at distance at least $r$ of the boundary of one of these subdivisions.
Then the union over all projections of the bad components over these overlapping subdivisions contains  $V\setminus V[K'r]$.

The overall conclusion is that we can  obtain good estimates on $V[r]$ by using a finite number of orthogonal projections and (possibly overlapping) subdivisions of the projection bases by cubes of size $Kr$.

\subsubsection{} Our main motivation for introducing this formalism is to obtain geometric information on wedge products of the form $v_n^{-1}[V_n]\wedge T^q$, where $T$ the Green current of an endomorphism of $\pk$  (or more generally any   closed positive current  wedgeable with $[V_n]$). We denote by $v_n^{-1}[V_n[r]]\wedge T^q$ the restriction of $v_n^{-1}[V_n]\wedge T^q$ to ${V_n[r]}$. What we  want is to estimate the proportion of the mass of $v_n^{-1}[V_n]\wedge T^q$ that is concentrated on the bounded geometry part $V_n[r]$, and ultimately show that 
in certain situations,  this proportion is close to 1 when $r$ is small, uniformly in $n$.

 This issue was addressed for curves in dimension 2 in various contexts \cite{isect, birat, dt-saddle, ddg2}. A crucial technical point in these papers is the validity of the  volume  estimate $$\unsur{v_n}\vol(V_n \setminus V_n[r]) = O(r^2).$$ To give some credit to the geometric assumptions that we will make in Section \ref{sec:transversex}, we prove this estimate in the general case in  \S \ref{subs:SA} below.
 Unfortunately, it leads to the desired result only when $q=1$. The details are given in \S \ref{subs:isect}.
It seems that for larger $q$, the --presumably optimal-- volume estimate that we obtain is not enough in itself to  control the mass of $v_n^{-1}[V_n[r]]\wedge T^q$, and that
 a finer understanding of the geometry of the ``bad part'' $V_n\setminus V_n[r]$ is required.

In \S \ref{subs:isect} we also show that for any $q$,  if $v_n^{-1}[V_n]$ converges to the current of integration over a $q$ dimensional analytic set, then  the expected control on the geometry of $v_n^{-1}[V_n]\wedge T^q$ is true.

\subsubsection{} Finally, let us quote a result that we will use several times.

\begin{thm}[Sibony-Wong \cite{sibony-wong}]\label{thm:sibonywong}
let $g$ be a holomorphic function defined in the neighborhood of the origin in $\cc^q$, which admits a holomorphic continuation to a neighborhood of $\bigcup_{L\in E} L\cap B(0,R)$, where $E\subset \pp^{q-1}$ is a set of lines through the origin,  of measure $\geq 1/2$
(relative to the Fubini-Study volume on $\pp^{q-1}$).

Then there exists a constant $C_{SW}>0$ such that $g$ extends to a holomorphic function on $B(0, C_{SW}R)$, and furthermore 
\begin{equation}\label{eq:sibonywong}
  \sup_{B(0, C_{SW}R)}\nolimits \abs{g} \leq
\sup_{\bigcup_{L\in E} L\cap B(0,R)}\nolimits \abs{g}.
\end{equation}
\end{thm}

\subsection{Volume estimates}\label{subs:SA}

In \cite{dinh-lamin}, Dinh proved the following theorem.

\begin{thm}[Dinh]\label{thm:dinh}
Let $q<k$ and $\iota_n: \pp^q\cv\pk$ be a sequence of holomorphic mappings, of generic degree 1.
Let $V_n= \iota_n(\pp^q)$,  $v_n$ be the volume of $V_n$, and $S_n= \unsur{v_n}[V_n]$.

Then every cluster value of the sequence of currents $(S_n)$ is  woven.
\end{thm}

Notice that if $f$ is an endomorphism of $\pk$ and $L$ is a generic line, then $V_n=f^n(L)$ satisfies the above assumptions.
As announced above, here we make this result more precise as follows.

\begin{thm}\label{thm:SA}
As in Theorem \ref{thm:dinh}, 
let $q<k$ and $\iota_n: \pp^q\cv\pk$ be a sequence of holomorphic mappings, of generic degree 1,
$V_n= \iota_n(\pp^q)$, $v_n$ be the volume of $V_n$, and   $S_n= \unsur{v_n}[V_n]$.

Denote by $V_n[r]$ the part of $V_n$ with bounded geometry at scale $r$,  and
$S_n[r]=\unsur{v_n}[V_n[r]]$ . Then there exists a constant $C$ such that
\begin{equation}\label{eq:SA}
\m(S_n-S_n[r])\leq Cr^2
\end{equation}
\end{thm}

\begin{proof}
The plan of the proof is close to that of \cite{dinh-lamin}, but we need to make things more explicit (see also \cite{dt-lamin}). We  assume that $q>1$; the case $q=1$ is more classical
and essentially contained in Lemma \ref{lem:moving_b}.
Let $C$ denote  a ``constant", that may vary from line to line,
independently of $n$ and $r$.


As explained in \S\ref{subsub:cube}, we will approximate $V_n[r]$ by a union of graphs over family of cubes over a family of projections $\pi_j$ in general position.  To construct these graphs, we will first construct graphs over generic 1-dimensional slices and glue those into $q$-dimensional graphs by using the Sibony-Wong Theorem \ref{thm:sibonywong}.

\medskip

Let us get into the details.
Fix  a linear subspace $I$ of dimension $(k-q-1)$, such that $I\cap V_n=\emptyset$ for all $n$, and let  $\pi:\pk\setminus I\cv \pp^q$ be the  projection of center $I$. Then $\pi\circ \iota_n:\pp^q\cv \pp^q$ is a holomorphic map of topological degree $v_n$. Equivalently, $\pi\rest{V_n}$ is a branched covering of
degree $v_n$.
Let also $E_n\subset V_n$ be the (Zariski closed) set of points $x$ such that $\#\iota_n^{-1}(\set{x})>1$.

In the projection base we work locally so we may assume that we are in  a bounded subset $\om$ of $\cc^q$,
equipped with  its standard metric.
If $L\subset \cc^q$ is  line, and $U\subset \cc^q$ is an open set , with $L\cap U\neq\emptyset$ we say that an {\em irreducible} component of $\pi^{-1}(L)\cap V_n$  is {\em good} if it is a graph $\Gamma$ over $L\cap U$ with $\int_\Gamma \omega_{\pp^k}\leq A$. Here $A$ is a constant whose value is chosen as follows:
we require that if $\Gamma$ is a good component  over $B(x,r)\subset L$ (with $r$ small enough, say, $r\leq D/10$, where $D$ is the constant of \S \ref{subsub:bddgeom}), then the diameter of $\Gamma \cap B(x, r/2)$ is less than $D/2$. This is possible because by the area-diameter lemma of  \cite{briend-duval2},
the diameter of $\Gamma \cap B(x, r/2)$ is bounded by a constant depending  only on $A$.

We let $g_n(L, U)$ be the number of good components  of $\pi^{-1}(L)\cap V_n$ over $L\cap U$, and $b_n(L, U) = v_n- g_n(L, U)$ the number of bad components counted with multiplicity.
  
We denote by $dz$ the Lebesgue measure on $L$. We have the following lemma (which will be proven later).
 
\begin{lem}\label{lem:moving_b}
There exists a constant $C$ such that if
 $L$ is a generic line and $r>0$, then
$$
\int_{L\cap \om} b_n(L, B(z,r)) dz \leq C v_nr^2.
$$
\end{lem}

We take $U$ to be a cube of size $r$.
Recall from \S\ref{subsub:cube} the notions of good and bad components over $U$.
We define $g_n(U)$ (resp. $b_n(U)$) to
be the number of good (resp.  bad) components $V_n$ over $U$ (resp. counting multiplicity). 

The following lemma relates $q$-dimensional good components  and good components over varying lines.
If $x\in \cc^q$ and $\delta\in \pp^{q-1}$ we let $L(x, \delta)$ be the line through $x$ with direction $\delta$. We simply denote by $d\delta$ the Fubini-Study volume element  on $\pp^{q-1}$.

\begin{lem}
There exists a constant $M$ depending only on the dimension such that if $Q$ is a cube of size $2r$, and $x\in Q$ is a.e. point
\begin{equation}\label{eq:SW1}
b_n(Q) \leq 2 \int_{\pp^{q-1}} b_n\lrpar{L(x,\delta), B(x,Mr)} d\delta.
\end{equation}
In particular
\begin{equation}\label{eq:SW2}
b_n(Q) \leq \frac{2}{r^{2q}}\int_{\pp^{q-1}\times Q} b_n\lrpar{L(x,\delta), B(x,Mr)} d\delta\hspace{.1em}  dx.
\end{equation}
\end{lem}

\begin{proof}[Proof of the lemma]  Let $M$ be greater than ${4\sqrt{q}}/{C_{SW}}$, where $C_{SW}$ is as in Theorem \ref{thm:sibonywong}. This is adjusted so that if $x\in Q$, $B(x, C_{SW}Mr)\supset Q$.

Let now $x\in Q$ be any point  belonging neither to the set of critical values of $\pi\rest{V_n}$ nor to $\pi(E_n)$. Then $\pi^{-1}\set{x}\cap V_n =  \set{y_1, \ldots , y_{v_n}}$ has cardinality $v_n$  and
in some neighborhood of $x$, $\pi\rest{V_n}$ admits $v_n$ inverse branches $g_j$, with $g_j(x)= y_j$.

Let $\varphi:\set{1, \ldots, v_n}\times \pp^{q-1} \cv \set{0,1}$ be defined by
$\varphi(j, \delta) = 0$ if  $\pi^{-1}(L(x,\delta))\cap V_n$ admits a good component over $B(x,Mr)\cap L(x,\delta)$ issued from $y_j$, and $\varphi(j, \delta) = 1$ otherwise.
Thus $b_n\lrpar{L(x,\delta), B(x,Mr)} = \sum_{j=1}^{v_n}\varphi(j, \delta)$.

For fixed $j$, by   Theorem \ref{thm:sibonywong}, if $\int_{\pp^{q-1}} \varphi(j,\delta) d\delta< 1/2$, then there is a good component of $V_n$ over $Q$ attached to $y_j$ --the control on the diameter comes from \eqref{eq:sibonywong}.
Therefore
$$b_n(Q) \leq \# \set{j,\ \int_{\pp^{q-1}} \varphi(j,\delta) d\delta> 1/2}
\leq 2 \sum_{j=1}^{v_n}\int_{\pp^{q-1}} \varphi(j,\delta) d\delta = 2\int_{\pp^{q-1}}
b_n\lrpar{L(x,\delta), B(x,Mr)},$$ which is the first part of the statement. The second part is obvious.
\end{proof}

We can now conclude the proof of the theorem.
As in \S\ref{subsub:cube}, consider a family of overlapping subdivisions  of a neighborhood of $\om$ in the projection base into cubes of size $2r$ (for convenience we put $K'=1$). If $\qq$ is one of them, and $V_{n, \qq}$ is the union of good components, we have
\begin{equation}\label{eq:subd}
\unsur{v_n} \bra{ [V_n]- [V_{n,\qq}], \pi^* \omega ^q_{\cc^q} } \leq C  {r^{2q}}\frac1{v_n}  {\sum_{Q\in \qq} b_n(Q)}.
\end{equation}
If we are able to show that the right hand side is a $O(r^2)$, then by taking the union of bad components relative to the overlapping subdivisions, for each projection $\pi_j$, and by using the fact that
$\sum_j \pi_j^* \omega ^q_{\pp^q} \geq c\omega_{\pp^k}^q$, we obtain the desired estimate on $\vol(V_n\setminus V_n[r])$.


Let $G(1,q)$ be the space of lines in $\cc^q$, endowed with its natural isometry-invariant  measure $\nu$. Let
$\widetilde G(1,q) = \set{(z,L),\ z\in L}$ be the tautological bundle over $G(1,q)$. It also possesses a natural measure, which, abusing slightly,  we denote by $dz\otimes \nu$, where $dz$ denotes Lebesgue measure on $L$. There is a natural diffeomorphism $\cc^q\times \pp^{q-1} \cv \widetilde G(1,q)$, which sends
$dx\hspace{.1em}d\delta$ to $dz\otimes \nu$ (up to a multiplicative constant).

Summing \eqref{eq:SW2} over all squares and changing variables
 we obtain the following estimate of the total number of bad components:
\begin{align}\label{eq:bnq}
\sum_{Q\in \qq} b_n(Q)  &\leq \frac{2}{r^{2q}} \sum_{Q\in \qq} \int_{\pp^{q-1}\times Q} b_n\lrpar{L(x,\delta), B(x,Mr)} d\delta\hspace{.1em}  dx \\ \notag
&=  \frac{2}{r^{2q}} \int_{\pp^{q-1}\times \om} b_n\lrpar{L(x,\delta), B(x,Mr)} d\delta\hspace{.1em}  dx\\ \notag
&=\frac{2}{r^{2q}} \int_{\set{(z,L)\in \widetilde G(1,q), \ z\in \om}}b_n\lrpar{L, B(z,Mr)} dz\hspace{.1em}d \nu(L)\\
&\leq C v_n r^{2-2q} \notag,
\end{align}
where the last inequality follows from Lemma \ref{lem:moving_b} and the  fact that the measure of the set of lines intersecting $\om$ is finite.  This, together  with \eqref{eq:subd}, completes the proof.
\end{proof}

\begin{proof}[Proof of Lemma \ref{lem:moving_b}]
We are considering $\pi:V_n\cap \pi^{-1}(L\cap \om)\cv L\cap\om$,  a branched covering of degree $v_n$
 between Riemann surfaces (recall that $I\cap V_n=\emptyset$, hence the projection is proper). In this situation, counting bad components is the same as counting critical values with multiplicity, plus discarding components of too large volume.

Again, we introduce subdivisions by cubes. In $\cc$, there exist 4 overlapping subdivisions $(\qq^i)_{i=1\ldots 4}$ by   squares of size $4r$ with the property that
 for any $z$, there exists $i(z)$ such that
$B(z,r)$ in contained in one square of $\qq^{i(z)}$.

We denote by $Q^i(z)$ the square of $\qq^i$ containing $z$ (there is an ambiguity for points at the boundaries of the subdivisions, but these have  zero Lebesgue measure).
Then $b_n(L, B(z,r))\leq b_n(L, Q^{i(z)}(z))\leq \sum_{i=1}^4  b_n(L, Q^{i }(z))$.
So $$\int_{L\cap \om} b_n(L, B(z,r)) dz \leq
\sum_{i=1}^4 \int_{L\cap \om}  b_n(L, Q^{i}(z)) = \sum_{i=1}^4 \sum_{{Q\in \qq^i }, {Q\cap\om\neq \emptyset } }(4r)^2 b_n(L, Q),$$ and we are left with proving that
 the number of bad components over $\qq^i$ is $O(v_n)$.

This is  classical; we recall the details for completeness.
Bad components are of two kinds: components with ramification points and graphs with too large volume.
To count ramified components, notice that since  $L\cap Q$ is simply connected, $b_n(L, Q)$ is not greater than the sum of the multiplicities of the critical points of
 $\pi\circ\iota_n\rest{(\pi\circ\iota_n)^{-1} (L)}$ on $(\pi\circ\iota_n)^{-1} (L\cap Q)$. Now  recall that $\pi\circ\iota_n$  is a holomorphic mapping $\pp^q\cv \pp^q$, of topological degree $v_n$, hence of degree $v_n^{1/q}$. Thus the degree of its critical set equals $(q+1)  (v_n^{1/q}-1)$ and the degree of the preimage of a line is $v_n^{(q-1)/q}$. We infer that the total number of critical points of  $\pi\circ\iota_n\rest{(\pi\circ\iota_n)^{-1} (L)}$, with multiplicity is not greater than  $(q+1)v_n$.

Regarding  bad components of the second kind, simply note that
since the volume of $V_n\cap \pi^{-1}(L)$ is $v_n$, there are not more than $\frac{v_n}{A}$ of them. This finishes the proof.
\end{proof}

\begin{rmk}
Another criterion for wovenness is given in  \cite[Theorem 5.6]{dinh-lamin}. The estimate \eqref{eq:SA}  also holds in this case.
\end{rmk}

\subsection{Geometric intersection} \label{subs:isect}
The following is a rather straightforward adaptation of \cite{isect}.

\begin{thm}\label{thm:isect}
Let $q=1$ and assume that  $V_n=\iota_n(\pu)$ and 
 $S_n$ are as in Theorem \ref{thm:dinh}. Let $T$ be a closed positive
current of bidegree (1,1) with continuous potential. Then $T\wedge \unsur{v_n}[V_n]$ is carried by the bounded geometry part of $V_n$, uniformly in $n$, that is, there exists a function $\e(r)$ tending to zero as $r\cv 0$ such that for all $n$,
\begin{equation}\label{eq:concentr}
\m(T\wedge  S_n[r])\geq 1-\e(r).
\end{equation}
\end{thm}

The condition on the potential of $T$ could actually be significantly relaxed,
along the lines of \cite{birat, ddg2}.

\begin{proof}
This is  identical to the 2 dimensional case, so we just outline the main steps of the argument (this will also be needed in the proof of Theorem \ref{thm:transversex}). The problem is local so it is enough to prove the mass estimate in a neighborhood some ball. Let $u$ be a potential of $T$ there. Fix $\e>0$.

As in the proof of Theorem \ref{thm:SA}, consider  linear projections $\pi_j$ to $\cc^q$ in general position, and generic subdivisions $\qq_j$ of the projection bases into cubes of size $r$. Then we get a subdivision $\mathcal{C}$ of $\cc^k$ by affine cubes of size $O(r)$, and a variety $V_{n,\mathcal{C}}$, which in each cube is a union of graphs relative to the $\pi_j$, satisfying the estimate $ \m(S_n-  S_{n,\mathcal{C}}) \leq Cr^2$
(with $S_{n,\mathcal{C}}=\unsur{v_n}[V_{n,\cC}]$). It is enough to prove that $\m(T\wedge (S_n-S_{n,\cC}))<\e$ when $r$ is small enough.

Recall from \S \ref{subsub:cube} the notation $\lambda C$ for the homothetic of $C$ of factor $\lambda$.  By \cite[Lemma 4.5]{isect} there exists $\lambda>0$ depending only on $\e$, and a translate of $\cC$ (still denoted by $\cC$), 
possibly depending on $n$, so that
$$\m\left(T\wedge S_n\rest{\bigcup_{ C\in\cC}C\setminus \lambda C}\right)
<\frac{\e}{2}.$$

Now, $\lambda$ being fixed, we estimate the mass in $\bigcup_{ C\in\cC}\lambda C$. For this, we let $\psi_\cC$ be a cutoff function, $0\leq \psi_\cC\leq 1$, equal to 1 in the neighborhood of every $\lambda C$, $C \in \cC$. It is possible to choose such a $\psi$ with $\norm{dd^c\psi}_{L^\infty}\leq O(r^{-2})$. Let $S_{n,C} = S_{n, \mathcal{C}}\rest{C}$.
In each cube we write
\begin{align}\label{eq:ipp}
 \bra{T \wedge  (S_n-S_{n,C}), \psi}& = \int_C \psi (dd^c u)  \wedge (S_n-S_{n,C})  \\ \notag
&=\int_C \left(u-u(\mathrm{center}(C)\right)  dd^c\psi  \wedge(S_n-S_{n,C}) \\ \notag &\leq
C \unsur{r^2} \omega(u,r)\m(S_n- S_{n,{C}}) ,
\end{align}
where $\omega(u,r)$ is the modulus of continuity of $u$. Then, summing over all cubes and using the volume estimate shows that
$$\m\left(T \wedge  (S_n - S_{n, \mathcal{C}})\rest{\bigcup_{ C\in\cC}
\lambda C}\right)\leq
C \omega(u,r).$$ Thus if $r$ is small enough this is less than $\frac{\e}{2}$ and we are done. Carefully inspecting the proof reveals that $\e$ depends only on $r$ (see \cite[Remark 4.7]{isect}
\end{proof}

Another instance where we are able to control $T^q\wedge [V_n]$ is the following one.

\begin{thm}\label{thm:isect_algebraic}
Let $q<k$, $V_n=\iota_n(\pp^q)$ and $S_n$ be as in Theorem \ref{thm:dinh}. Let $T$ be a closed positive
current of bidegree (1,1) with continous potential and assume that $S_n$ converges to the current of integration over a $q$-dimensional analytic cycle. Then $T^q\wedge S_n$ is carried by the bounded geometry part of $V_n$, uniformly in $n$, i.e. \eqref{eq:concentr} holds.
\end{thm}

\begin{proof}
Let $V = \sum \alpha_j V_j$ be such that $[V] = \lim \unsur{v_n}[V_n]$. The probability measure $T^q\wedge [V]$ gives no mass to proper analytic subsets of $V$ so it is carried by the regular part $\mathrm{Reg}( V)$.
Let $U\subset \mathrm{Reg}( V)$ be a ball, and $U'$,
 be a     tubular neighborhoods of $U$ in $\pk$,  so small that
$U'\cap V=U$. Fix $\e>0$. We will show that if $r$ is small and $n$ is large, then  
\begin{equation}\label{eq:cover}
\m\left(T^q\wedge \unsur{v_n}\left[V_n[r] \cap{U'}\right]\right)\geq  \m(T^q\wedge[V\cap{U'}])-\e.\end{equation} 
A simple covering argument then leads to \eqref{eq:concentr}.

\medskip

As in the proof of Theorem \ref{thm:SA}, consider a linear subspace $I$ of dimension $(k-q-1)$, such that $I\cap V$, as well as $I\cap V_n$ are empty for all $n$, and let  $\pi:\pk\setminus I\cv \pp^q$ be the  projection of center $I$.
We may further assume that the fibers of the projection are transverse to $V$ in $U$, so that locally we   view $\pi$ as a projection onto $U$. Let $Q(x,r)$ be the cube of center $x$ and radius $r$.

 By a slight variation on  \cite[Lemma 5.2]{dinh-lamin} which we  explain below, for any $\eta>0$ there exists a proper analytic subset $C_\eta$  such that if $x\notin C_\eta$ and $r$ is small enough 
(uniformly on compact subsets of $U\setminus C_\eta$), then 
$V_n$ contains at least $v_n(1-\eta)$ good components over $Q(x,2r)$. Note that the fiber $\pi^{-1}(x)$ possibly intersects $V$ in several other points. Write 
$V_n\rest{\pi^{-1}( Q(x,r))} = G_n+ B_n$, where $G_n$  denotes  the union of strong good  components (for convenience, from now on we drop the ``strong''). Let $\alpha$ be the coefficient of the component containing $U$ in the cycle $V$, so that $[V]\rest{U'} = [V\cap U']= \alpha[U]$.

Recall that the diameter $D$ of good components can be chosen arbitrarily small (provided $r$ is).
Consider another  tubular neighborhood $U''\subset U'$ of $U$, with $U''\cap V=U$, and  choose $D$ so small that 
 any good component intersecting $U''$ is contained in $U'$.

We know that ${v_n^{-1}}[V_n\cap{U'}]\cv \alpha [U]$. We want to show that asymptotically, at least $(\alpha-\eta)v_n$ good components over $Q(x,r)$ are contained in $U'$. For this, we count good components by projecting their volume, i.e. integrating $\pi^*\omega^q$ on them,  and write
\begin{align*}
 \bra{\unsur{v_n}[G_n\cap{U''}], \pi^*\omega^q}&=
\bra{\unsur{v_n}[V_n\cap{U''}\cap \pi^{-1}(Q(x,r))], \pi^*\omega^q} - \bra{\unsur{v_n}[B_n\cap{U''}], \pi^*\omega^q}\\
&\geq \bra{\unsur{v_n}[V_n\cap{U''}
\cap \pi^{-1}(Q(x,r))], \pi^*\omega^q} - \eta\int_{Q(x,r)} \omega^q \\
&\underset{n\cv\infty}\longrightarrow (\alpha-\eta) \int_{Q(x,r)} \omega^q .
\end{align*}
This means that  
asymptotically, at least $(\alpha-\eta)v_n$ good components over $Q(x,r)$ intersect $U''$.
Thanks to our choice of the diameter $D$
we conclude that for large $n$ there is a current $[G'_n]$
made of at least $(\alpha-\eta)v_n$ good components over $Q(x,r)$,
 entirely contained in $U'$, and such that $\liminf(v_n^{-1}[{G}'_n])\geq
 (\alpha-\eta)[Q(x,r)]$. 
 
 Since $T$ has continuous potential we infer that
 for large  $n$,  $\m(v_n^{-1}[{G}'_n]\wedge T^q)\geq (\alpha -2\eta)\m([Q(x,r)]
 \wedge T^q)$. As explained in \S \ref{subsub:cube},
  $G'_n$ is contained in $V_n[r]$ (again for convenience we put $K'=1$), so this may be rephrased as
$$\m\left(\unsur{v_n} [V_n[r] \cap \pi^{-1}(Q(x,r))\cap U']\wedge T^q\right)\geq (\alpha -2\eta)\m([Q(x,r)] \wedge T^q).$$ Finally, since $[U]\wedge T^q$ carries no mass on $C_\eta$, we can adjust the constant $\eta$ and
 cover a set of large  $([U]\wedge T^q)$-mass with finitely many disjoint cubes of radius $r$ avoiding $C_\eta$,
and  conclude that \eqref{eq:cover} holds.

\medskip

It just remains to explain how  \cite[Lemma 5.2]{dinh-lamin} should be modified to construct the  analytic subset $C_\eta$. Following Dinh's notation, there exists two positive closed currents $S_\infty$ and $\Omega_\infty$  of bidimension $(q,q)$, and a constant $\nu(\eta)$ 
such that if the mass of $S_\infty$ (resp. $\Omega_\infty$)  in $B(x, 2r)$  is less than
$r^{2-2q} \nu$, then $V_n$ admits $v_n(1-\eta)$ good components over $Q(x,r)$, as desired. We see that this is true for small $r$ as soon as $\nu(S_\infty, x)< \nu_0 = :\nu/(2\sqrt{2q})^{2-2q}$  
(resp. $\nu(\Omega_\infty, x)< \nu_0$). Here $\nu$ denotes the Lelong number; for the value of $\nu_0$, note  that  
$Q(x,r)\subset B(x,\sqrt{2q}r)$. Therefore it is enough to put
$C_\eta = \set{x, \  \nu(S_\infty, x)\geq \nu_0 \text{ or }  
\nu(\Omega_\infty, x)\geq \nu_0}$,  
which is an analytic set by Siu's Theorem \cite{siu}.
\end{proof}

\section{Lower estimates for transverse expansion} \label{sec:transversex}

The main theorem in this section is the following.
 Recall that if $V$ is a subvariety in $\pk$, we denote by $V[r]$ the part of $V$ with bounded geometry at scale $r$. 

\begin{thm}\label{thm:transversex}
Let $f$ be a holomorphic endomorphism of $\mathbb{P}^k$, and $T$ be its Green current.
Assume that the following holds: 

\begin{narrower}
$(H_q)$  For a.e. linear subspace $L$ of dimension $q$, there exists a subexponentially decreasing sequence $(r_n)$
such that 
\begin{equation}\label{eq:unifgeom}
 \m\left(\unsur{d^{nq}}T^q\wedge \big[(f^nL)[r_n]\big]\right) \underset{n\cv\infty}\longrightarrow 1.
\end{equation}
\end{narrower}

Then for $\sigma_{T^q}$-a.e. $x$ and (Lebesgue) a.e
$q$-dimensional complex linear subspace $ V \subset T_x\pk$, if
$v\in V$ is a non-zero vector
then
\begin{equation}\label{eq:transversex}
\limsup_{n\cv\infty}\unsur{n}\log\norm{df^n_x(v)}\geq \frac{\log d}{2}.
\end{equation}
\end{thm}

By subexponentially decreasing, of course we mean that  $\lim \frac{\log{r_n}}{n} = 0$. In dimension 2, this result is related to the work of  De Thélin \cite{dt-saddle}, which itself
relies on techniques introduced by the author for studying the dynamics of birational mappings \cite{birat}. We also borrow some arguments from Dinh and Sibony \cite{ds-pl}. We believe that $(H_q)$  always holds. This opinion is of course supported by
the analysis  of the previous section.

\begin{rmk}\label{rmk:transversex}
The following facts are consequences of the proof of the theorem. We leave the reader fill the details.
\begin{enumerate}
 \item Under assumption $(H_q)$ our proof actually 
  shows the more precise result  that for every $\e>0$ there exists a set of integers $\nn_\e\subset \nn$ of density at least $1-\e$ such that 
$$
\liminf_{\nn_\e\ni n\cv\infty}\unsur{n}\log\norm{df^n_x(v)}\geq \frac{\log d}{2}.
$$ In this case $\nn_\e$ depends on $(x,V)$.  Likewise,  to get \eqref{eq:transversex}, it is enough to require \eqref{eq:unifgeom} along  a subsequence $(n_j)$.   
\item If we only assume that $\liminf \frac{\log{r_n}}{n} \geq -\alpha$, for some $0\leq \alpha< \frac{\log d}{2}$, then
we obtain a similar conclusion, with the right hand side of \eqref{eq:transversex} replaced by  $\frac{\log d}{2}-\alpha$. This variation is sufficient to imply Corollary \ref{cor:fatoumajor}.
\item The result can be localized as follows: if there exists an open set $U\subset \pk$ which is a union of $q$ dimensional linear spaces, such that $(H_q)$ holds for $L\subset U$, then \eqref{eq:transversex} holds in $U$.
\end{enumerate}
\end{rmk}

\begin{cor}\label{cor:fatoumajor}
If $(H_q)$ holds, then for $\sigma_{T^q}$ a.e. $x\in J_q\setminus J_{q+1}$, the Fatou subspace has dimension $k-q$ at $x$. Moreover,  if $v$ is any  tangent vector at $x$, not belonging to $\mathcal{F}_x$, then $v$ satisfies the expansion property \eqref{eq:transversex}. In other words, Conjecture \ref{conj:fatoumain} holds for $q$. 

Additionally, $T^q$ is decomposable a.e. on $J_q\setminus J_{q+1}$ and for $\sigma_{T^q}$ a.e. $x\in  J_q\setminus J_{q+1}$,  $\mathcal{T}^q_x = \mathcal{F}_x$.
\end{cor}

\begin{proof}
We already know by Corollary \ref{cor:fatouminor} that the dimension of the Fatou subspace is at least $k-q$. Assume that the inequality is strict on a set of positive trace measure $E$. Then if $x\in E$ and $V\subset T_x\pk$ is a $q$ dimensional subspace, $\dim (V\cap\mathcal{F}_x) \geq 1$. This of course contradicts Theorem \ref{thm:transversex}. 

From Corollary \ref{cor:fatouminor} again, we know that  $\mathcal{T}^q_x\subset \mathcal{F}_x$. In particular  it has dimension $k-q$ and  $T^q$ is decomposable a.e.

Part {\em ii.} of Conjecture \ref{conj:fatoumain} is then immediate. By genericity, there exists a supplementary subspace $U$ to $\mathcal{F}_x$ satisfying \eqref{eq:transversex}. Now if $v$ is any vector outside $\mathcal{F}_x$, it admits a decomposition as $v = e+u$, with $e\in \mathcal{F}_x$ and $u\in U\setminus\set{0}$. Applying $df^n_x$ and using the fact that $\norm{df^n_x(e)}$ grows subexponentially then gives the result.
\end{proof}

It was shown in Theorem  \ref{thm:isect} that $(H_1)$  holds for any sequence $r_n\cv0$. Therefore we have:

\begin{cor}\label{cor:q=1}
Conjecture \ref{conj:fatoumain} is  true for $q=1$.
\end{cor}

The following corollary applies for instance when there exists a $q$-dimensional algebraic attractor.  The proof will be given afterwards.

\begin{cor}\label{cor:attractor}
Assume that there exists an open set $U\subset \pk$, which is a union of $q$-dimensional subspaces, with the property that there exists a $q$-dimensional subvariety $V$ of $\pk$  such that every cluster value of $(f^n)_*(\sigma_{T^q}\rest{U})$ is concentrated on $V$. Then the conclusions of 
Theorem \ref{thm:transversex} hold in $U$. 
\end{cor}

\begin{proof}[Proof of  Theorem \ref{thm:transversex}]~

\noindent{\bf Step 1.} If $V$ is a subspace of $T_x\pk$, we put
$$ \chi(x,V) = \limsup_{n\cv\infty} \inf_{{v\in V}, \ { \norm{v}=1}} \unsur{n}\log\norm{df^n_x(v)}.$$
Let us first prove  that the conclusion of the theorem is true if the following holds:
\begin{equation}\label{eq:expansionL}
\text{for a.e. }q\text{-dimensional subspace } L, \text{ for }
T^q\wedge [L]\text{ a.e. }x,  \ \chi(x,T_xL)\geq \frac{\log d}{2}.
\end{equation}
Of course $\sigma_{T^q}$ is an average of $T^q\wedge [L]$. Still, it is not obvious to deduce  \eqref{eq:transversex} from \eqref{eq:expansionL}, for
 in \eqref{eq:expansionL} the direction  of expansion depends on the measure.

\medskip

It is no loss of generality to assume
that we work in  a ball  $B\subset \cc^k$ endowed with its usual metric (we still denote the K\"ahler form by $\omega$).
Let ${G(q,k)}$ be the space of $q$-dimensional subspaces of $\cc^k$, and $\nu$ be its Haar measure.  If $V\in G(q,k)$,
let $\omega_V$ be the current defined by
$\omega_V = \int_{V^\perp} [V+u] d\leb(u)$ (integration over the family of subspaces of direction $V$). Notice that $\int \omega_V d\nu(V) = \omega^{k-q}$. Now suppose  that the conclusion of the theorem is false. Then there exists a measurable set
$E\subset B\times G(q,k)$ with positive $\sigma_{T^q}\otimes \nu$ measure
such that  $\chi(x,V)< \log{d}/2$ when
 $(x,V)\in E$. Let $E_V = \set{x, \ (x,V)\in E}$; this is a set of points where the expansion property fails in direction $V$. There exists $\alpha>0$ and a set $A_\alpha$ of positive $\nu$ measure  such that if $V\in A_\alpha$,
  $\sigma_{T^q}(E_V)\geq \alpha$. We may further assume that for $V\in A_\alpha$, generic subspaces $L\subset B$ of direction $V$ satisfy   \eqref{eq:expansionL}.

If $S$ is a set and $N\geq 2$, we let  $S^{[N]} = S^N/\mathfrak{S}_{N}$ be the set of subsets of $S$ with cardinality $N$.
Notice that  $G(q,k)^{[N]}$ is  endowed with a  natural measure $\nu_{N}$ derived from $\nu$.
We know from  linear algebra  that there  exists $N_1$ depending only on $q$ and $k$ such  that if $V_1, \ldots , V_{N_1}$ is a collection
 of subspaces in general position, then
$\sum_{i=1}^{N_1}\omega_{V_i}>0$, i.e. $\sum_{i=1}^{N_1}\omega_{V_i} \geq \e\omega$ for some $\e>0$.
The set of such $\set{V_1, \ldots,V_{N_1}}$ is open
and of full $\nu_{N_1}$ measure in $G(q,k)^{[N_1]}$. Furthermore, if $N_2>N_1$, the set of collections $\mathcal{V} = \set{V_1, \ldots, V_{N_2}}$ of subspaces of cardinality $N_2$  such that for every
$\set{V_{j_1}, \ldots, V_{j_{N_1}}}\subset \mathcal{V} $,
$\sum_{i=1}^{N_1}\omega_{V_{j_i}}>0$, is also open and of full $\nu_{N_2}$ measure.

We conclude that for every $N_2>N_1$, there exists a collection $\mathcal{V}  = \set{V_1, \ldots , V_{N_2}}\subset A_\alpha$ with the property that for every
$\set{V_{j_1}, \ldots, V_{j_{N_1}}}\subset \mathcal{V}$,
$\sum_{i=1}^{N_1}\omega_{V_{j_i}}>0$.

Fix $N_2>N_1/\alpha$ and $\mathcal{V}$ as above. Since $\int \sum_{j=1}^{N_2} \mathbf{1}_{E_{V_j}} \sigma_{T^q} \geq \alpha N_2> N_1$ and $\m(\sigma_{T^q})\leq 1$, we infer that there exists a set of positive trace measure of points belonging to at least $N_1$ subsets $E_{V_j}$. Since $\mathcal{V}^{[N_1]}$ is finite,
there exists a particular collection $\set{V_{j_1}, \ldots, V_{j_{N_1}}}$ such that $F= \bigcap_{i=1}^{N_1} E_{V_{j_i}}$ has positive trace measure. Finally, since $\sum_{i=1}^{N_1} \omega_{V_{j_i}}\geq \e\omega$ we conclude that there exists $1\leq i\leq N_1$ with the property that  $(T^q\wedge \omega_{V_{j_i}})(E_{V_{j_i}})\geq (T^q\wedge \omega_{V_{j_i}})(F)>0$.
Now, $T^q\wedge\omega_{V_{j_i}}$ is an average of measures $T^q\wedge [V_{j_i}+u]$ for which by assumption, \eqref{eq:expansionL} holds a.e. for vectors belonging to $V_{j_i}$. This contradicts the definition of $E_{V_{j_i}}$.

\medskip

\noindent{\bf Step 2.} Proof of \eqref{eq:expansionL} for generic $L$.

If $L$ is a  $q$-dimensional subspace, then by Bézout's Theorem $f^n(L)$ has degree (hence volume) $d^{nq}$. If furthermore $L$ is generic, then $f^n\rest{L}: L\cv f^n(L)$ is  1-1 outside some subvariety (a birational map). Indeed this happens when for some $x\in L$,  $L$ meets $f^{-n}(\set{f^n(x)})$ only at $x$, which  clearly holds outside a proper Zariski closed set. We take such a $L$ and  assume that it satisfies \eqref{eq:unifgeom}.

The plan of the proof is the following.
 Let $V_n = f^n(L)$. As usual, we realize the bounded geometry part of $V_n$ as a union of graphs over subdivisions by cubes. We  introduce a family of dynamically defined bad components and by  reconsidering the proof of Theorem \ref{thm:isect}, we check that  discarding them
  does not affect the mass estimate \eqref{eq:unifgeom}.
Finally, for the remaining part of $V_n$ we obtain  good expansion estimates leading to \eqref{eq:expansionL}.

Fix $\e>0$ and an integer $n$.  The estimate in \eqref{eq:unifgeom} is local so we work in a ball.
Fix projections $\pi_j$, subdivisions of the projection bases  $\mathcal{Q}_j$, and the resulting subdivision $\cC$ as in the proof of Theorem \ref{thm:isect}, except that the size of the cubes is now $r_n$. Consider the family of homothetic cubes $\lambda C$.
If the $\mathcal{Q}_j$ are well positioned,
there exists $\lambda<1$ {\em depending only on $\e$} such that 
$(T \wedge S_n)(\cC\setminus \la \cC)<\e/2$. This value of $\la$ is fixed from now on. For each projection $\pi_j$, we form the variety $V_{n, j,\mathcal{Q}_j}$ made  of the good components of $V_n$ over $\mathcal{Q}_j$, and  let $S_{n,j,\mathcal{Q}_j} = \unsur{d^{nq}}[V_{n, j,\mathcal{Q}_j}]$, as usual.

If $\Gamma$ is such a good component,  abusing slightly we  denote by
$f^{-n}(\Gamma)\cap L$  the proper transform of $\Gamma$ under $(f^n\rest{L})^{-1}$. Then
$f^n: f^{-n}(\Gamma)\cap L\cv L$ is a biholomorphism.   Indeed,
$f^{-n}(\Gamma)\cap L$ and $\Gamma$ are smooth, and $f^n: f^{-n}(\Gamma)\cap L\cv \Gamma$
is both finite and birational.
It is well-known that it must be a biholomorphism in this case. Indeed the critical set, if nonempty, is a hypersurface, and the local structure of $h$ near a smooth point of this hypersurface is that of $(x_1, \ldots, x_k)\mapsto (x_1^\alpha, \ldots, x_k)$ for some $\alpha\geq 2$. 

From now on the  inverse of $f^n\rest{f^{-n}(\Gamma)\cap L}$ will be  denoted by $f_{-n}$.

If $\Gamma$ is a component of $V_{n, j,\mathcal{Q}_j}$, consider the integral
$$I_n(\Gamma) = \int_\Gamma (f_{-n})^*\omega\wedge\o^{q-1} = \int_{f_{-n}(\Gamma)} \o\wedge (f^n)^* \o^{q-1}.$$
 Since the $f_{-n}(\Gamma)$ are disjoint open subsets of $L$, we infer that
$$\sum_{\Gamma \text{ comp. of }V_{n, j,\mathcal{Q}_j}}I_n(\Gamma)  \leq \int_L\o\wedge (f^n)^* \o^{q-1} = d^{n(q-1)}.$$ Therefore,
$$\# \set{\Gamma\text{ component of }V_{n, j,\mathcal{Q}_j}, \ I_n(\Gamma) \geq \frac{1}{d^n}}\leq {d^{nq}}.$$
Discard these components from $V_{n, j,\mathcal{Q}_j}$, and let $V'_{n, j,\mathcal{Q}_j}$ (resp  $S'_{n, j,\mathcal{Q}_j}$) be the remaining variety  (resp. current). Since we have removed at most $d^{nq}$ graphs, we have that
\begin{equation}
\label{eq:r_n}
\bra{S_{n, j,\mathcal{Q}_j} - S'_{n, j,\mathcal{Q}_j}, \pi_j^*\omega_{\pp^q}^q}\leq C r_n^{2q}.
\end{equation}

We can now form the currents $S_{n,\cC}$ (resp. $S'_{n,\cC}$), by taking,
in each cube $C\in \cC$, the union of the
components of $S'_{n, j,\mathcal{Q}_j}$ (resp.$S'_{n, j,\mathcal{Q}_j}$).
As explained in \S\ref{subsub:cube}, \eqref{eq:unifgeom} implies that $\m(T^q\wedge S_{n,\cC})\cv1$ as $n\cv\infty$ (again, for convenience we put $K=1$).

From \eqref{eq:r_n} we infer that $\m(S_{n,\cC}-S'_{n,\cC})\leq C r_n^{2q}$. We can now estimate $$\m\left((T^q\wedge
S_{n,\cC}- T^q\wedge S'_{n,\cC})\rest{\bigcup_{C\in \cC}\la C}\right)$$ by applying exactly the same reasoning that in Theorem \ref{thm:isect}, except that \eqref{eq:ipp} is replaced by a sequence of $q$ integration by parts, leading to the inequality
$$\m\left((T^q\wedge
S_{n,\cC}- T^q\wedge S'_{n,\cC})\rest{\bigcup_{C\in \cC}\la C}\right) \leq C \unsur{r_n^{2q}} \omega(u,r_n)^q \m(S_{n,\cC}-S'_{n,\cC})\leq C \omega(u,r_n)^q.$$
We conclude that when $n$ is large enough, $\m(T^q\wedge S'_{n,\cC})\geq 1-\e$.

\medskip

We now construct a set $A_{n, \e}$, with $([L]\wedge T^q)(A_{n,\e})>1-\e$, such that if $x\in A_{n, \e}$ and $v$ is  a unit vector tangent to $L$, then $\norm{df^n_x(v)}\geq C(\e) d^{(1-\e)n/2}$.
For this, observe that $(f^n)_*([L]\wedge T^q) = S_n\wedge T^q$.
Therefore,    $(f_{-n})_*((S'_{n, \cC}\wedge T^q )\rest{\la \cC})$ is
 a measure dominated by $[L]\wedge T^q$, with mass larger than $ 1-\e$. Let then
 $A_{n,\e}  =\bigcup_{\Delta \text{ comp. of }S'_{n,\cC}}  f_{-n}(\Delta \cap \la C)$.

The control on the derivative comes from the following lemma.

\begin{lem}\label{lem:derivative}
Let $\Delta$ be a component of $S_{n,\cC}$. Then there exists $C(\e)$ such that
$\norm{df_{-n}}\leq C(\e)\frac{d^{-n/2}}{r_n}$ on $\Delta\cap \la C$.
\end{lem}

This result being assumed for the moment, we can finish the proof. Recall first that $\frac{\log r_n}{n}\cv 0$,
so ${r_n}d^{n/2}\geq d^{(1-\e)n/2}$ for large $n$. Next, if we set $$B_{\e} = \set{x, \ \exists \nn_{\e} \text{ of density } \geq 1-\sqrt{\e}, \ \text{s.t. } \forall n\in \nn_\e, \ x\in A_{n,\e}},$$
it is an exercise  (see \cite[Lemma 6.5]{ddg3}) to  show that $([L]\wedge T^q)(B_\e) \geq 1-\sqrt{\e}$.
Thus, $[L]\wedge T^q$-a.e. point belongs to  $B_\e$ for some $\e$, and we are done.
\end{proof}

Lemma \ref{lem:derivative} will itself follow from a result of independent interest.

\begin{prop}\label{prop:sibony_wong}
Let $D \subset \cc^q$ be a bounded convex domain in $\cc^q$,  endowed with its natural metric and
K\"ahler form $\beta$, and let  $(X, \o_X)$ be a compact Hermitian  manifold (where $\o_X$ denotes the $(1,1)$ form associated to the metric).

Then  for every compact subset $K\subset D$, there exists a constant $C(K)$
 such that for every holomorphic mapping $h:D\cv X$ we have
\begin{equation}\label{eq:sw}
\lrpar{\mathrm{diam} (h(K))}^2\leq C \int_D h^*\omega_X\wedge \beta^{q-1}.
\end{equation}
\end{prop}

This, combined with the Cauchy estimates in coordinate charts, gives estimates on the derivative of $h$. For $q=1$ this is  a rough version of
 the area-diameter inequality of Briend and Duval \cite{briend-duval2}, and for $q>1$ this is merely a reformulation of ideas due to Dinh and Sibony \cite{ds-pl}. We include the proof for convenience.

\begin{proof}[Proof of Proposition \ref{prop:sibony_wong}]
Since $X$ is compact it suffices to prove the result when the integral on the right hand side of \eqref{eq:sw} is small enough.

Let $d\theta$ be the  Fubini-Study  volume element on the space $\set{L_{x,\theta}, \theta\in \mathbb{P}^{q-1}}$ of complex lines through each $x\in D$. Let $\alpha_x$ be the current defined by $\alpha_x = \int[L_{x,\theta}] d\theta$ and observe that $\beta^{q-1}  = \int\alpha_x dx$ (up to a normalization factor which we assume is 1).

We say that an affine line $L_{x,\theta}$ is  $A$-good if $\int h^*\omega_X \wedge [L_{x,\theta}] \leq A
\int_D h^*\omega_X\wedge \beta^{q-1}$, where $A$ is a constant to be fixed later.
Likewise we say that $x\in D$ is $A$-good if among all lines through $x$,
the measure of the set of $A$-good ones is larger than $1/2$.
There exists a universal constant $C_1$ such that if $A \geq  \frac{C_1}{r^{2q}}$, each ball of radius $r$ contains an $A$-good point (argue by contradiction).

\medskip

 Let now $R>0$ be smaller than $\dist(K,\fr D)/4$. Define
$M>0$ to be the infimum  of the moduli of the annuli   $ (L\cap D)\setminus (L\cap B(y,R))$, where $B(y,R)$ is any ball of radius $R$ intersecting $K$, and $L$ is any line through $y$. Let $r= \frac{C_{SW}R}{2}$, where $C_{SW}$ is the constant appearing in the Sibony-Wong Theorem \ref{thm:sibonywong}, and $A =
\frac{C_1}{r^{2q}}$ as above.

Cover $K$ with a finite family of balls $B(x,r)$. The required number of course depends only on $K$.
Each of these balls
contains an $A$-good point $y$, and $B(y, C_{SW}R)\supset B(x,r)$.
Now, by the Briend-Duval area-diameter estimate, for every line through $y$  we have
$$\mathrm{diam}(h(L\cap B(y,R)))^2\leq \frac{\mathrm{Area}
(h(L\cap D))}{\mathrm{mod}(L\cap B(y,R),L\cap D)} \leq
\unsur{M}\int_L h^*\omega_X\leq \frac{A}{M}\int_D h^*\omega_X\wedge \beta^{q-1}.$$
Thus, if $\lrpar{\frac{A}{M}\int_D h^*\omega_X\wedge \beta^{q-1}}^{1/2}$ is less than the diameter of a coordinate chart of $X$,  Theorem \ref{thm:sibonywong} applies, and in particular
we obtain that
$h\rest{B(y,C_{SW}R)}$ takes its values in the chart, with the same estimate on the diameter. Since these balls cover $K$ the proof is finished.
\end{proof}

\begin{proof}[Proof of Lemma \ref{lem:derivative}]
By scaling, we may assume that the cube has size 1. This affects the derivative
by a factor $\frac{1}{r_n}$.

If $\Delta$ is a component of $S_{n,\cC}$, it is the restriction to a cube $C$ of a component of some
$S_{n,j,\mathcal{Q}_j}$, that is, a graph $\Gamma$ of a function $\gamma$ over a cube $Q\in \mathcal{Q}_j$ in $\cc^q$, and satisfying $I_n(\Gamma)\leq \unsur{d^n}$. Notice that $\pi_j\circ\gamma = \mathrm{id}$.  It will be enough to estimate the derivative of $f_{-n}$ on $\Gamma\cap \pi_j^{-1}(\la Q)$.

Write now $f_{-n}\rest{\Gamma} = (f_{-n}\circ\gamma)\circ\pi_j$. The  derivative of $\pi_j$ is
 uniformly bounded.
To deal with  that of $(f_{-n}\circ\gamma)\rest{\lambda Q}$,
we use Proposition \ref{prop:sibony_wong}.
By assumption on $ \Gamma$, we have that
$$I_n(\Gamma) = \int_{\Gamma} (f_{-n})^*\omega_L\wedge \omega_{ \Gamma}^{q-1}\leq \unsur{d^{n}}.$$
 Recall that we were working in a ball of $\cc^k$ so that,   we can freely consider
 $\omega_\Gamma$ as being the restriction to $\Gamma$ of the natural K\"ahler form in $\cc^k$.
Now, since $\Gamma$ is a graph, we have that $\omega_{\Gamma}\geq   \pi_j^* \omega_Q$, or equivalently
$\gamma^*\omega_\Gamma\geq \omega_Q$, so we infer that $$\int_{Q}(f_{-n}\circ \gamma)^*\omega_L\wedge\omega_Q^{q-1} \leq
\int_{Q}(f_{-n}\circ \gamma)^*\omega_L\wedge (\gamma^*\omega_{  \Gamma})^{q-1} =
\int_{\Gamma} (f_{-n})^*\omega_L\wedge \omega_{\Gamma}^{q-1} \leq
 d^{-n}.$$
Consequently from Proposition \ref{prop:sibony_wong} we conclude that
 the diameter of $(f_{-n}\circ\gamma)(\la Q)$ is bounded by
$C(\lambda)d^{-n/2}$. The constant depends only on the scaling factor $\lambda$, hence ultimately on $\e$.

Finally, working in charts and using the Cauchy estimates, we conclude that
$\norm{df_{-n}\rest{{\Gamma}\cap \la C}}\leq C(\e)d^{-n/2}$, which was the desired result.
\end{proof}

\begin{proof}[Proof of Corollary \ref{cor:attractor}]
We need to show that for a generic $q$-dimensional linear subspace $L\subset U$, \eqref{eq:unifgeom} holds. In the open subset of the Grassmanian $G(q,k)$ consisting of subspaces contained in $U$, consider a smooth probabiliy measure $m$, and the associated current $\Sigma_0 = \int [L] dm(L)$. Since $T^q\wedge \Sigma_0 \ll\sigma_{T^q}$ we infer that any cluster value of $(f^n)_* (T^q\wedge \Sigma_0)$ is concentrated on $V$.

Let $\Sigma_n =\unsur{d^{nq}}(f^n)_*\Sigma_0$. Since $(f^n)_*(T^q\wedge \Sigma_0) = T^q\wedge \Sigma_n$ and   $T$ has continuous potential,  the cluster values of $(f^n)_* (T^q\wedge \Sigma_0)$ are of the form  $T^q\wedge \Sigma_\infty$, with $\Sigma_\infty$a cluster value  of $\Sigma_n$.

The following lemma is certainly well-known.

\begin{lem}\label{lem:charge}
If $S$ is a closed positive current on $\pk$ of bidimension $(q,q)$ that gives no mass to a  complete pluripolar set $P$, then neither does $T^q\wedge S$.
\end{lem}

Let $\mathcal{C}_V$ be the cone  of currents of integration on  cycles supported on $V$.
Let $\Sigma_\infty$ be as above and decompose  $\Sigma_\infty$ as $\Sigma_\infty= \Sigma_V+ \Sigma_\infty'$, where $\Sigma_V\in\mathcal{C}_V$ and
$\Sigma_\infty'$  gives no mass to $V$. By the above lemma $T^q\wedge \Sigma_\infty'=0$ which by Bézout's Theorem  implies that $\Sigma_\infty'=0$. Thus we conclude that $\Sigma_\infty$ is a current of integration supported on $V$.

Now recall that $\mathcal{C}_V$ is an extremal face of the cone of positive closed currents, in the sense that if $\Sigma\in \mathcal{C}_V$ and $S\leq \Sigma$, then $S\in \mathcal{C}_V$ \cite{lelong-extremal}. Since  $\Sigma_n = \int\frac{[f^n(L)]}{d^{nq}} dm(L)$  is the barycenter of a measured family of positive closed currents, converging to $\mathcal{C}_V$, it is an exercise to show that for $m$-a.e. $L$, $\unsur{d^{nq}}[f^n(L)]$ converges to $\mathcal{C}_V$ as well.

From this and Theorem \ref{thm:isect_algebraic} we conclude that $(H_q)$ holds for a.e. $L\subset U$, and the proof is complete.
\end{proof}

\begin{proof}[Proof of Lemma \ref{lem:charge}](compare \cite[Prop. 1.2]{ddg2})
This is a local problem, so we work in a ball $B$.
Write $P = \set{\psi=-\infty}$ for some negative psh function $\psi$. By assumption $\sigma_S(P)=0$.
By replacing $\psi$ with $\chi\circ\psi$, where $\chi$ is a slowly growing
convex increasing function with $\lim_{-\infty}\chi = -\infty$, we can actually assume that $\psi\in L^1(\sigma_S)$. The following version of the Chern-Levine-Nirenberg inequality is true (see e.g. \cite[Thm A.3.2]{ds-survey}):
if $K\subset B$ is a  relatively compact open set, and the $(u_i)_{i=1}^q$ are bounded psh functions in $B$, then
$$\m_K(\psi dd^cu_1\wedge \cdots \wedge dd^cu_q\wedge S) \leq C(K) \norm{\psi}_{L^1(\sigma_S)} \norm{u_1}_{L^{\infty}(B)}\cdots \norm{u_q}_{L^{\infty}(B)}.$$ It follows that $\psi\in L^1_{\rm loc}(T^q\wedge S)$, and the result follows.
\end{proof}

\end{document}